\documentclass[onefignum,onetabnum]{siamart220329}

\usepackage{lipsum}
\usepackage{amsfonts}
\usepackage{graphicx}
\usepackage{epstopdf}

\usepackage{algorithm}
\usepackage{algorithmicx}
\usepackage[noend]{algpseudocode}

\algnewcommand\algorithmicclass{\textbf{class}}
\algdef{SE}[FUNCTION]{Class}{EndClass}%
   [2]{\algorithmicclass\ \textproc{#1}\ifthenelse{\equal{#2}{}}{}{(#2)}}%
   {\algorithmicend\ \algorithmicclass}%

\makeatletter
\ifthenelse{\equal{(\vec{A})LG@noend}{t}}%
  {\algtext*{EndClass}}
  {}%
\makeatother

\ifpdf
  \DeclareGraphicsExtensions{.eps,.pdf,.png,.jpg}
\else
  \DeclareGraphicsExtensions{.eps}
\fi

\newsiamremark{remark}{Remark}
\newsiamremark{example}{Example}
\newsiamremark{hypothesis}{Hypothesis}
\crefname{hypothesis}{Hypothesis}{Hypotheses}
\newsiamthm{claim}{Claim}

\headers{Optimal rational matrix function approximation}{T. Chen, A. Greenbaum, C. Musco, and C. Musco}

\title{Low-memory Krylov subspace methods for optimal rational matrix function approximation%
\thanks{
\textbf{Funding:} This material is based on work supported by the National Science Foundation under Grant Nos. DGE-1762114, CCF-2045590, and CCF-2046235 and by an Adobe Research grant.}}

\author{Tyler Chen\thanks{New York University, \href{mailto:tyler.chen@nyu.edu}{\ms{tyler.chen@nyu.edu}}}
\and Anne Greenbaum\thanks{University of Washington, \href{mailto:greenbau_uw.edu}{\ms{greenbau@uw.edu}}}
\and Cameron Musco\thanks{University of Massachusetts Amherst, \href{mailto:cmusco@cs.umass.edu}{\ms{cmusco@cs.umass.edu}}} 
\and \linebreak[4] Christopher Musco\thanks{New York University, \href{mailto:cmusco@nyu.edu}{\ms{cmusco@nyu.edu}}}
}

\usepackage{amsopn}

\crefname{section}{Section}{Sections}
\crefname{ineq}{}{}
\creflabelformat{ineq}{#2{\upshape(#1)}#3}

\usepackage{marginnote}

\usepackage{enumitem}

\usepackage{tikz}
\usetikzlibrary{calc}
\usepackage{graphicx}
\usepackage[caption=false]{subfig}

\usepackage{booktabs}
\usepackage{array}

\usepackage{afterpage}
\usepackage{pdflscape}

\newcommand{\Tyler}[1]{{\color{red}Tyler: #1}}

\DeclareMathOperator*{\argmin}{argmin}

\usepackage{bbm}

\usepackage{mathtools}

\newcommand{\ms}{\texttt}
\usepackage{bm}
\usepackage{dsfont}
\newcommand{\bOne}{\ensuremath{\mathds{1}}}

\newcommand{\R}{\mathbb{R}}

\renewcommand{\d}[1]{\ensuremath{\mathrm{d}#1}}

\renewcommand{\vec}{\mathbf}

\newcommand{\T}{\textup{\textsf{T}}}
\newcommand{\cT}{\textup{\textsf{T}}}
\newcommand{\ii}{\bm{i}}

\usepackage{xspace}

\newcommand{\lan}{\textup{\textsf{lan-FA}}}

\newcommand{\lanopt}{\textup{\textsf{lan-OR}}}

\allowdisplaybreaks

\usepackage[normalem]{ulem}

\ifpdf
\hypersetup{
  pdftitle={Low-memory Krylov subspace methods for optimal rational matrix function approximation},
  pdfauthor={T. Chen, A. Greenbaum, C. Musco, and C. Musco}
}
\fi

\begin{document}

\maketitle

\begin{abstract}
We describe a Lanczos-based algorithm for approximating the product of a rational matrix function with a vector. 
This algorithm, which we call the Lanczos method for optimal rational matrix function approximation (Lanczos-OR), returns the optimal approximation from a given Krylov subspace in a norm depending on the rational function's denominator, and can be computed using the information from a slightly larger Krylov subspace.
We also provide a low-memory implementation which only requires storing a number of vectors proportional to the denominator degree of the rational function. 
Finally, we show that Lanczos-OR can be used to derive algorithms for computing other matrix functions, including the matrix sign function and quadrature based rational function approximations.
In many cases, it improves on the approximation quality of prior approaches, including the standard Lanczos method, with little additional computational overhead.
\end{abstract}



\section{Introduction}
Krylov subspace methods (KSMs) are among the most powerful algorithms for computing approximations to $f(\vec{A})\vec{b}$ when $\vec{A}$ is an $n\times n$ real symmetric matrix and $f:\R\to\R$ is an arbitrary function.
Such methods construct an approximation to $f(\vec{A})\vec{b}$ that lies in the Krylov subspace 
\begin{equation*}
\mathcal{K}_k
:= \operatorname{span}\{\vec{b}, \vec{A}\vec{b}, \ldots, \vec{A}^{k-1}\vec{b} \}, 
\end{equation*}
and only need to access $\vec{A}$ through matrix-vector products.
This means KSMs are well suited for large-scale computations where storing $\vec{A}$ in fast memory is infeasible.

In the special case that $f(x) = 1/(x-z)$ for some $z\in\mathbb{C}$, KSMs such as conjugate gradient (CG) \cite{hestenes_stiefel_52}, minimum residual (MINRES) \cite{paige_saunders_75}, and quasi-minimum residual (QMR) \cite{freund_92} are able to provide \emph{optimal} approximations to $f(\vec{A})\vec{b}$ from $\mathcal{K}_k$ while storing just a few vectors of length $n$.
For general functions $f$, however, the situation is murkier. 
General purpose KSMs, like the well-known Lanczos method for matrix function approximation (Lanczos-FA) \cite{saad_92}, are not known to return an optimal or near-optimal approximation to $f(\vec{A})\vec{b}$ from $\mathcal{K}_k$, except for a few cases like the exponential \cite{druskin_greenbaum_knizhnerman_98}. 
In fact, even work on weaker spectrum dependent bounds remains somewhat ad-hoc, including for the basic case of rational functions \cite{eshof_frommer_lippert_schilling_van_der_vorst_02,ilic_turner_simpson_09,frommer_guttel_schweitzer_14,frommer_guttel_schweitzer_14a,chen_greenbaum_musco_musco_22}.

Moreover, in terms of computational cost, to return an approximation to $f(\vec{A})\vec{b}$ from $\mathcal{K}_k$, methods for general functions like Lanczos-FA either (i) store $k$ vectors of length $n$, or 
(ii) store a constant number of vectors of length $n$, but increase the number of matrix-vector products by a factor of two \cite{boricci_00a,frommer_simoncini_08}. Lower memory methods have been studied for specific classes of functions \cite{guttel_schweitzer_21}. For instance, for Stieltjes or analytic functions, restarting methods are a potential alternative to saving all of the $k$ vectors generated by Lanczos \cite{afanasjew_eiermann_ernst_guttel_08,niehoff_06,frommer_guttel_schweitzer_14a,frommer_guttel_schweitzer_14,ilic_turner_simpson_09}.
However, restarting can discard useful information from the Krylov subspace, possibly delaying convergence.

In this paper, we address the above issues with existing KSMs by describing an optimal algorithm with good memory performance for the important case of rational functions. 
We call the method the Lanczos method for optimal rational matrix function approximation (Lanczos-OR). Our method applies to \emph{any} rational function $f$. 
If the degrees of the numerator and denominator of $f$ are each at most $d$ (typically a small constant), then Lanczos-OR produces {optimal} approximations to $f(\vec{A})\vec b$ (in a certain norm\footnote{As dicussed in the next section, we prove optimality in a norm that depends on the rational function being approximated, but this norm is closely related to e.g., the more standard $2$-norm or $\vec{A}$-norm. Lanczos-OR performs well experimentally for these norms as well.}) from the span of $\mathcal{K}_k$, using at most $k + d$ matrix-vector products. In the special case when the denominator matrix is positive definite, Lanczos-OR is equivalent to the optimal Galerkin projection method from \cite[Section 4]{lopez_simoncini_06} and if $f(x) = 1/(x-z)$, the CG, MINRES, and QMR iterates are obtained as special cases.

Prior work in \cite{lopez_simoncini_06} largely viewed Lanczos-OR as a method of theoretical interest, that could possibly help explain the behavior of Lanczos-FA. 
In contrast, we argue that Lanczos-OR is a useful algorithm in and of itself, by showing how its iterates can be computed efficiently.
In addition to only requiring $d$ more matrix-vector products than the standard Lanczos-FA method, we provide an implementation of Lanczos-OR that requires storing just $2d + 4$ vectors of length $n$. Therefore, for a fixed rational function, the storage costs do not grow with the iteration $k$. 
Our approach can also be used for computing the Lanczos-FA approximations to rational matrix functions, avoiding storage costs growing with $k$ in that widely used method.

Beyond rational functions, we show that Lanczos-OR can be used to derive algorithms for approximating other functions. In particular, we derive ``induced'' rational approximations via integral representations of functions like the matrix sign function.
While not provably optimal, these induced Lanczos-OR approximations tend to perform well in practice.
In fact, on problems where Lanczos-FA exhibits erratic behavior, the Lanczos-OR induced approximations tend to have nicer behavior.

\subsection{The Lanczos algorithm and some basic Krylov subspace methods}
\label{sec:background}

KSMs for symmetric matrices  are often based on the Lanczos algorithm.
Given a symmetric matrix $\vec{A}$ and vector $\vec{b}$, the Lanczos algorithm run for $k$ iterations constructs an orthonormal basis \( \vec{Q} := [ \vec{q}_0 , \ldots , \vec{q}_{k-1} ] \) such that the first $j\leq k$ columns form a basis for the Krylov subspace
\begin{equation*}
\mathcal{K}_j= \operatorname{span}\{ \vec{b}, \vec{A}\vec{b}, \ldots, \vec{A}^{j-1}\vec{b} \}
= \{ p(\vec{A}) \vec{b} : \deg(p) < j \}.
\end{equation*}
Moreover, the basis vectors satisfy a symmetric three term recurrence,
$$
\vec{A} \vec{Q} = \vec{Q} \vec{T} + \beta_{k-1} \vec{q}_{k} \vec{e}_{k-1}^\T.
$$
Here $\vec{e}_{k-1}$ is the standard basis vector with a one in the last entry and $\vec{T}$ is symmetric tridiagonal with  diagonals $(\alpha_0, \ldots, \alpha_{k-1})$ and off diagonals $(\beta_0, \ldots, \beta_{k-2})$  which are also computed by the algorithm. 

In our analysis it will be useful to consider the recurrence that would be obtained if the Lanczos algorithm were run to completion. 
In exact arithmetic, for some $K\leq n$, $\beta_{K-1} = 0$ in which case the algorithm terminates. 
Then the final basis \( \widehat{\vec{Q}} := [ \vec{q}_0 , \ldots , \vec{q}_{K-1} ]\) and symmetric tridiagonal $\widehat{\vec{T}}$ with  diagonals $(\alpha_0, \ldots, \alpha_{K-1})$ and off diagonals $(\beta_0, \ldots, \beta_{K-2})$
satisfy
 a three-term recurrence 
\begin{equation*}
    \vec{A} \widehat{\vec{Q}} = \widehat{\vec{Q}} \widehat{\vec{T}}.
\end{equation*}
Since the columns of $\widehat{\vec{Q}}$ are orthonormal, we have that $   \widehat{\vec{T}}
    = \widehat{\vec{Q}}^\T \vec{A} \widehat{\vec{Q}}$,
from which we easily see that, after any number of iterations $k$, $\vec{T} = \vec{Q}^\T \vec{A} \vec{Q}$. Note that $\vec{Q} = [\widehat{\vec{Q}}]_{:,:k}$ and $\vec{T} =[\widehat{\vec{T}}]_{:k,:k}$.
Note also that, for any shift $z\in\mathbb{C}$, $    (\vec{A} - z\vec{I}) \widehat{\vec{Q}} = \widehat{\vec{Q}} (\widehat{\vec{T}} - z\vec{I})$.
In other words, the Krylov subspaces generated by $(\vec{A},\vec{b})$ and $(\vec{A}-z\vec{I},\vec{b})$ coincide, and the associated tridiagonal matrices are easily related by a diagonal shift.

\subsection{Notation}
We denote the complex conjugate of $z$ by $\overline{z}$.
Matrices and vectors are denoted by bold upper and lowercase letters, respectively.
We use \emph{zero-indexed numpy style slicing} to indicate entries. 
Specifically, \( [\vec{B}]_{r:r',c:c'} \) denotes the submatrix of \(\vec{B}\) consisting of rows \( r \) through \( r'-1 \) and columns \( c \) through \( c'-1 \). 
If any of these indices are equal to \( 0 \) or \( n \), they may be omitted and if \( r' = r+1  \) or \( c' = c+1 \), then we will simply write \( r \) or \( c \).
For example, \( [\vec{B}]_{:,:2} \) denotes the first two columns of \( \vec{B} \) (corresponding to indices 0 and 1), and \( [\vec{B}]_{3,:} \) denotes the fourth row (corresponding to index 3).
Throughout, $\vec{A}$ will be a real symmetric matrix.
We denote the set of eigenvalues of $\vec{A}$ by $\Lambda$ and define $\mathcal{I} := [\lambda_{\textup{min}}, \lambda_{\textup{max}}]$.
Without loss of generality, we assume that $\|\vec{b}\|_2 = 1$, where $\vec{b}$ is the vector to which the rational function is applied. 

\section{Optimal rational function approximation}

We now describe an optimal iterate for approximating $r(\vec{A})\vec{b}$ when $r(x)$ is a rational function whose denominator is nonzero at the eigenvalues $\Lambda$ of $\vec{A}$.
We will describe a low-memory implementation of this algorithm in \cref{sec:low_mem} that can also be used to efficiently compute Lanczos-FA approximations to $r(\vec{A})\vec{b}$.

\begin{definition}
\label{def:lanczos_or}
Let \( r(x) \) be a rational function written as \( r(x) = M(x)/N(x) \), where \( N(x) \) is a polynomial with leading coefficient one and \( M(x) \) is a polynomial sharing no common factors with \( N(x) \).
For any polynomial $R(x)$, define \( \tilde{M}(x) = M(x)R(x) \) and \( \tilde{N}(x) = N(x) R(x) \).
Then the Lanczos-OR iterate is defined as
\begin{equation*}
    \lanopt_k(r,R)
    :=  \vec{Q}  ([\tilde{N}(\widehat{\vec{T}})]_{:k,:k})^{-1} [\tilde{M}(\widehat{\vec{T}})]_{:k,:k} \vec{e}_0.
 \end{equation*}
\end{definition}

Those familiar with CG, MINRES, and the version of QMR for shifted Hermitian systems will note that these optimal algorithms are each obtained as special cases of Lanczos-OR.
Specifically, when $\vec{A}$ is positive definite, CG is obtained with \( r(x) = 1/x \) and \( R(x) = 1 \), MINRES is obtained with $r(x) = 1/x$ and \( R(x) = x \), and QMR is obtained if $r(x) = 1/(x-z)$ and \( R(x) = (x-\overline{z}) \).
In fact, we prove a more general optimality result for Lanczos-OR:
\begin{theorem}
\label{thm:lanczos_opt}
Given a rational function $r(x) = M(x)/N(x)$ as in \cref{def:lanczos_or}, choose a polynomial $R(x)$ so that $\vec{H} = \tilde{N}(\vec{A}) = N(\vec{A})R(\vec{A})$ is positive definite.
Then $\lanopt_k(r,R)$ is the $\vec{H}$-norm optimal approximation to $r(\vec{A})\vec{b}$ from $\mathcal{K}_k$; i.e.,
\begin{align*}
    \|r(\vec{A})\vec{b}- \lanopt_k(r,R) \|_{\vec{H}} = \min_{\vec{x}\in\mathcal{K}_k}
    \| r(\vec{A}) \vec{b} - \vec{x} \|_{\vec{H}}.
 \end{align*}
\end{theorem}

\begin{proof}
Since it lies in $\mathcal{K}_k$, the minimizer $\vec{x}_*$ of $\min_{\vec{x}\in\mathcal{K}_k} \| r(\vec{A}) \vec{b} - \vec{x} \|_{\vec{H}}$  can be written as $\vec{Q} \vec{c}_*$ for
\begin{equation*}
    \vec{c}_* = \argmin_{\vec{c}\in\R^k} \| r(\vec{A})\vec{b} - \vec{Q} \vec{c} \|_{\vec{H}}
    =
    \argmin_{\vec{c}\in\R^k} \| \vec{H}^{1/2} r(\vec{A})\vec{b} - \vec{H}^{1/2} \vec{Q} \vec{c} \|_2.
\end{equation*}
This is a standard least squares problem which yields
\begin{align}
    \label{eqn:OR_c}
    \vec{x}_* 
    = \vec{Q}\vec{c}_* 
    = \vec{Q}( \vec{Q}^\cT \vec{H} \vec{Q})^{-1} \vec{Q}^\cT \vec{H} r(\vec{A}) \vec{b}.
\end{align}
By definition, $\vec{H} = N(\vec{A})R(\vec{A}) $ so $\vec{H}r(\vec{A}) = M(\vec{A})R(\vec{A}) = \tilde{M}(\vec{A})$. Thus,
\begin{equation}
    \label{eqn:OR_QHrAb}
    \vec{Q}^\cT \vec{H} r(\vec{A}) \vec{b}
    = \vec{Q}^\cT \tilde{M}(\vec{A}) \vec{b}
    = \vec{Q}^\cT \widehat{\vec{Q}}\tilde{M}(\widehat{\vec{T}}) \widehat{\vec{Q}}^\cT \vec{b}
    = [\tilde{M}(\widehat{\vec{T}})]_{:k,:k} \vec{e}_0.
\end{equation}
Next, since $\vec Q$ consists of the first $k$ columns of $\widehat{\vec{Q}}$, 
\begin{equation}
    \label{eqn:QHQ}
    \vec{Q}^\cT \vec{H} \vec{Q}
    = \vec{Q}^\cT \tilde{N}(\vec{A}) \vec{Q}
    = [\widehat{\vec{Q}}^\cT \tilde{N}( \vec{A} ) \widehat{\vec{Q}}]_{:k,:k}
    =  [\tilde{N}(\widehat{\vec{T}}) ]_{:k,:k}.
\end{equation}
Plugging \cref{eqn:OR_QHrAb,eqn:QHQ} into \cref{eqn:OR_c}, and using that $\vec{x}_* = \vec{Q}\vec{c}_*$ yields the result.
\end{proof}

Even though \cref{thm:lanczos_opt} establishes that Lanczos-OR returns an optimal approximation to $r(\vec{A})\vec{b}$ in a non-standard norm, the $\vec{H}$-norm, this optimality already implies a number of nice properties. For example, it immediately implies that, up to a multiplicative factor independent of the iteration $k$, the Lanczos-OR iterates are comparable to the optimal 2-norm approximations to $r(\vec{A})\vec{b}$. Formally, we have: 
\begin{corollary}
\label{thm:2normopt}
Given a rational function $r(x) = M(x)/N(x)$ as in \cref{def:lanczos_or}, choose a polynomial $R(x)$ so that $\vec{H} = \tilde{N}(\vec{A}) = N(\vec{A})R(\vec{A})$ is positive definite.
Then,
\begin{equation*}
    {\|r(\vec{A})\vec{b}- \lanopt_k(r,R) \|_2}/{\|\vec{b}\|_2}
    \leq 
    \sqrt{\kappa(\vec{H})} \min_{\vec{x}\in\mathcal{K}_k} {\| r(\vec{A}) \vec{b} - \vec{x} \|_2}/{\|\vec{b}\|_2}.
\end{equation*}
\end{corollary}
\begin{proof}
Using basic properties of the $\vec{H}$-norm we have
\begin{align*}
    \|r(\vec{A})\vec{b}- \lanopt_k(r,R) \|_2
    &\leq \sqrt{\|\vec{H}^{-1}\|_2}
    \|r(\vec{A})\vec{b}- \lanopt_k(r,R) \|_{\vec{H}}
    \\&=\sqrt{\|\vec{H}^{-1}\|_2}\min_{\vec{x}\in\mathcal{K}_k}
    \| r(\vec{A}) \vec{b} - \vec{x} \|_{\vec{H}}
    \\&\leq \sqrt{\kappa(\vec{H})} \min_{\vec{x}\in\mathcal{K}_k} \| r(\vec{A}) \vec{b} - \vec{x} \|_2.
\end{align*}
\end{proof}

In \cref{sec:2normopt} we provide an experiment which suggests that the factor $\sqrt{\kappa(\vec{H})}$ may  be very pessimistic in some cases.

Based on \cref{thm:lanczos_opt}, we also obtain an a priori error bound involving the best scalar polynomial approximation to $r$ on the eigenvalues of $\vec{A}$, analogous to the well known minimax bounds for CG, MINRES, and QMR \cite{greenbaum_97} and \cite[Proposition 4.2]{lopez_simoncini_06}.
\begin{theorem}
Given a rational function $r(x) = M(x)/N(x)$ as in \cref{def:lanczos_or}, choose a polynomial $R(x)$ so that $\vec{H} = \tilde{N}(\vec{A}) = N(\vec{A})R(\vec{A})$ is positive definite.
Then, 
\begin{equation*}
    {\| r(\vec{A}) \vec{b} - \lanopt_k(r,R) \|_{\vec{H}}}/{\| \vec{b} \|_{\vec{H}}}
    \leq \min_{\deg(p) < k} \max_{\lambda \in \Lambda} | r(\lambda) - p(\lambda) |.
\end{equation*}
\end{theorem}

\begin{proof}
Since $\lanopt_k(r,R)$ is the $\vec{H}$-norm optimal approximation over the Krylov subspace, we have
\begin{equation*}
    \| r(\vec{A})\vec{b} - \lanopt_k(r,R) \|_{\vec{H}}
    = \min_{\vec{x}\in \mathcal{K}_k} \| r(\vec{A})\vec{b} - \vec{x} \|_{\vec{H}}
    = \min_{\deg(p)<k} \| r(\vec{A})\vec{b} - p(\vec{A})\vec{b} \|_{\vec{H}}.
\end{equation*}
Next, using the fact that $\vec{A}$ and $\vec{H}^{1/2}$ commute, we note that,
\begin{equation*}
    \| r(\vec{A})\vec{b} - p(\vec{A})\vec{b} \|_{\vec{H}}
    = \| (r(\vec{A})- p(\vec{A}))\vec{H}^{1/2}\vec{b} \|_2
    \leq \| r(\vec{A})- p(\vec{A}) \|_2 \| \vec{b} \|_{\vec{H}}.
\end{equation*}
Finally, the result follows from using the definition of the spectral norm to write
\begin{align*}
    \| r(\vec{A})- p(\vec{A}) \|_2
    = \| (r - p)(\vec{A}) \|_2
    = \max_{\lambda \in \Lambda} | r(\lambda)-p(\lambda) |.
\end{align*}
\end{proof}

\subsection{Efficient computation of the optimal iterate}
While its optimality and the resulting bounds above imply that the Lanczos-OR iterate should be a natural choice for rational function approximation, preferred over e.g. the standard Lanczos-FA approximation, it is not yet apparent that the Lanczos-OR iterate can be computed efficiently using a small number of matrix-vector multiplications. Naively, the iterate involves the terms  $[\tilde{N}(\widehat{\vec{T}})]_{:k,:k}$ and $[\tilde{M}(\widehat{\vec{T}})]_{:k,:k}\vec{e}_0$, and computing $\widehat{\vec{T}}$ requires running Lanczos to completion.
Fortunately these quantities can be computed efficiently.

\begin{lemma}
\label{thm:poly_That}
Suppose $p$ is a polynomial with $q := \deg(p) > 0$ and put $k':=k+\lfloor q/2 \rfloor$.
Then
\begin{equation*}
    [p(\widehat{\vec{T}})]_{:k,:k} = [p([\widehat{\vec{T}}]_{:k',:k'})]_{:k,:k}.
\end{equation*}
Moreover, $[p(\widehat{\vec{T}})]_{:k,:k}$ can be computed using the coefficients generated by $k+\lfloor (q-1)/2 \rfloor$ iterations of Lanczos.
\end{lemma}

\begin{proof}
It suffices to consider the case $p(x) = x^q$.
Let $\widehat{\vec{I}}_\ell$ be a $K\times \ell$ matrix whose top $\ell$ rows are an $\ell\times \ell$ identity and bottom $K-\ell$ rows are all zero. We have that
\begin{equation*}
    \widehat{\vec{T}} \widehat{\vec{I}}_\ell
    = [\widehat{\vec{T}}]_{:,:\ell}
    = \widehat{\vec{I}}_{\ell+1} [\widehat{\vec{T}}]_{:\ell+1,:\ell}.
\end{equation*}
Repeatedly applying this relation,
\begin{equation*}
    \widehat{\vec{T}}^j \widehat{\vec{I}}_k
    = 
    \widehat{\vec{I}}_{k+j}
    [\widehat{\vec{T}}]_{:k+j,k+j-1}
    \cdots
    [\widehat{\vec{T}}]_{:k+2,k+1}
    [\widehat{\vec{T}}]_{:k+1,k}
    = \widehat{\vec{I}}_{k+j} \vec{B}(k+j,k)
\end{equation*}
where we have defined
\begin{equation*}
    \vec{B}(k+j,k) := [\widehat{\vec{T}}]_{:k+j,k+j-1}
    \cdots
    [\widehat{\vec{T}}]_{:k+2,k+1}
    [\widehat{\vec{T}}]_{:k+1,k}.
\end{equation*}
Therefore, since $(\widehat{\vec{I}}_{k+j})^\cT \widehat{\vec{I}}_{k+j}$ is the $(k+j)\times (k+j)$ identity,
\begin{equation*}
    [\widehat{\vec{T}}^{2j}]_{:k,:k}
    = \vec{B}(k+j,k)^\cT \vec{B}(k+j,k)
\end{equation*}
and, since $(\widehat{\vec{I}}_{k+j-1})^\cT \widehat{\vec{T}} \widehat{\vec{I}}_{k+j-1} = [\vec{T}]_{:k+j-1,:k+j-1}$,
\begin{equation*}
    [\widehat{\vec{T}}^{2j-1}]_{:k,:k}
    = \vec{B}(k+j-1,k)^\cT  [\widehat{\vec{T}}]_{:k+j-1,:k+j-1} \vec{B}(k+j-1,k).
\end{equation*}
The expressions for $[\widehat{\vec{T}}^{2j}]_{:k,:k}$ and $[\widehat{\vec{T}}^{2j-1}]_{:k,:k}$ both depend only on $[\widehat{\vec{T}}]_{:k+j,:k+j-1}$, which can be obtained using $k+j-1$ matrix-vector products.
The first claim of the lemma follows by noting that $\lfloor (2j-1)/2 \rfloor = j-1$. To complete the lemma, note that for $\ell+1\leq k'$,
\begin{equation*}
    [\widehat{\vec{T}}]_{:k',:k'} \widetilde{\vec{I}}_\ell
    = [\widehat{\vec{T}}]_{:k',:\ell}
    = \widetilde{\vec{I}}_{\ell+1} [\widehat{\vec{T}}]_{:\ell+1,:\ell}
\end{equation*}
where $\widetilde{\vec{I}}_\ell$ is defined like $\widehat{\vec{I}}_\ell$ but is $k' \times \ell$.
The same argument as above then gives
\begin{equation*}
    ([\widehat{\vec{T}}]_{:k',:k'})^j \widetilde{\vec{I}}_k
    =\widetilde{\vec{I}}_{k+j}\vec{B}(k+j,k)
\end{equation*}
provided that \( k+j \leq k' \).
We therefore have that 
\(
    [\widehat{\vec{T}}^{q}]_{:k,:k} = [([\widehat{\vec{T}}]_{:k',:k'})^{q}]_{:k,:k}. 
\)
\end{proof}

We can therefore bound the number of matrix-vector products required to compute the Lanczos-OR iterates.
\begin{corollary}
Given a rational function $r(x) = M(x)/N(x)$ as in \cref{def:lanczos_or}, the Lanczos-OR iterate $\lanopt_k(r,R)$ can be computed using $\max\{\deg(\tilde{M})+1,k+\lfloor \deg(\tilde{N})/2\rfloor\}$ matrix-vector products, where $\tilde{M}$ and $\tilde{N}$ are as in \cref{def:lanczos_or}.
\end{corollary}

\begin{proof}
It is well-known that $[\tilde{M}(\widehat{\vec{T}})]_{:k,:k}\vec{e}_0 = \tilde{M}([\widehat{\vec{T}}]_{:k,:k})\vec{e}_0$ for any $k \geq \deg(\tilde{M})$ \cite{druskin_knizhnerman_89,saad_92}.
Thus, $[\tilde{M}(\widehat{\vec{T}})]_{:k,:k}\vec{e}_0$ can be computed using $\deg(\tilde{M})+1$ matrix-vector products.
Then, using \cref{thm:poly_That} we have that $[\tilde{N}(\widehat{\vec{T}})]_{:k,:k}$ can be computed using $k+\lfloor \deg(\tilde{N})/2\rfloor$ matrix-vector products.
\end{proof}

Since we assume that $N(x)$ is nonzero on the spectrum of $\vec{A}$, a simple way to ensure $\vec{H}$ is positive definite is to take \( R(x) = N(x) \) so that $\vec{H} = N(\vec{A})^2$.
However, in some situations, we may be able to use a lower degree choice for \( R(x) \), often resulting in a better conditioned $\vec{H}$.
For instance, in the case of symmetric linear systems, while one can always use MINRES ($r(x)=1/x$, $R(x) =x$), if \( \vec{A} \) is positive definite, then one typically would use CG ($r(x)=1/x$, $R(x) =1$). 
A simple way to obtain a lower degree choice
of $R(x)$ is to take only the terms in $N(x)$ which are indefinite.
\begin{definition}
    \label{def:Rstar}
    Given a rational function $r(x) = M(x)/N(x)$ as in \cref{def:lanczos_or}, factor \( N(x) \) as 
    \begin{equation*}
        N(x) = (x-z_1) \cdots (x-z_{d_1}) (x-{z_1'})(x-\overline{z_1'}) \cdots (x-{z_{d_2}'})(x-\overline{z_{d_2}'})
    \end{equation*}
    where $z_i\neq \overline{z}_j$ for all $i,j = 1,\ldots, d_1$ with $j \neq i$.
    Then \( R^*(x) \) is defined by
    \begin{equation*}
        R^*(x) = \xi (x-\overline{z_1})^{\alpha_1} \cdots (x-\overline{z_{d_1}})^{\alpha_{d_1}}
    \end{equation*} 
    where, for \( i=1,\ldots,d_1 \), $\alpha_i = 0$ if $z_i\in \mathbb{R} \setminus \mathcal{I}$  and $\alpha_i = 1$ otherwise and $\xi\in\{\pm 1\}$ is chosen so that $R^*(\lambda_{\textup{min}})N(\lambda_{\textup{min}}) > 0$.  
\end{definition}

\begin{lemma}
    Given a rational function $r(x) = M(x)/N(x)$ as in \cref{def:lanczos_or}, choose \( R^* \) as in \cref{def:Rstar}. Then \( \vec{H} = N(\vec{A}) R^*(\vec{A}) \) is positive definite.
\end{lemma}

\begin{proof}
    For each $i=1,\ldots, d_2$, $(x-{z_i'})(x-\overline{z_i'}) \geq 0$ for all $x \in \mathbb{R}$.  For each $z_i \in \mathbb{R} \setminus \mathcal{I}$, $i=1,\ldots, d_1$, $(x-z_i)$ does not change signs over $\mathcal{I}$.
The choice of $\xi$ ensures that $\tilde{N} (x) = N(x) R^{*}(x)$ is nonnegative throughout $\mathcal{I}$, and since, by assumption,
$N( \lambda ) \neq 0$ for any $\lambda \in \Lambda$, it follows that $\tilde{N} ( \lambda )$ is positive and therefore
that $\vec{H} = \tilde{N}(\vec{A})$ is positive definite.
\end{proof}

\section{Algorithms for other matrix functions}

In this section we discuss how Lanczos-OR can be used to derive algorithms for non-rational matrix functions, using the matrix sign function as a running example. We focus on the value of rational functions obtained from integral representations of a target function $f$. Such representations have been used in a range of past work on Krylov subspace methods to provide error bounds or estimates, and even to derive more advanced approximation schemes, such as restarted Lanczos-FA \cite{ilic_turner_simpson_09,frommer_guttel_schweitzer_14,frommer_guttel_schweitzer_14a,frommer_schweitzer_15,chen_greenbaum_musco_musco_22}. 

\subsection{Leveraging integral representations}
\label{sec:rational_func}
One possible use of Lanczos-OR is to approximate a rational matrix function $r(\vec{A})\vec{b}$, which is iitself an approximation to some non-rational matrix function $f(\vec{A})\vec{v}$.
For any output $\mathsf{alg}(r)$ meant to approximate $r(\vec{A})\vec{b}$, we have the following bound:
\begin{align}
\| f(\vec{A}) \vec{b} - \mathsf{alg}(r) \|
&\leq 
\nonumber
\| f(\vec{A}) \vec{b} - r(\vec{A}) \vec{b} \| +  \| r(\vec{A}) \vec{b} - \mathsf{alg}(r) \|
\\& \leq \underbrace{\|\vec{b}\|  \max_{\lambda \in \mathcal{I}} | f(\lambda) - r(\lambda) |}_{\text{approximation error}} 
    + \underbrace{\vphantom{\max_{\lambda \in \mathcal{I}}}\| r(\vec{A}) \vec{b} - \mathsf{alg}(r) \|}_{\text{application error}}.
    \label{eqn:triangle_ineq}
\end{align}
In many cases, very good or even optimal scalar rational function approximations to a given function on a single interval are known or can be easily computed. 
Thus, the approximation error term can typically be made small with a rational function of relatively low degree. At the same time, the bound is only meaningful if the approximation error term is small relative to the application error.

Rational function approximations commonly are obtained by discretizing an integral representation using a numerical quadrature approximation. 
For instance, the matrix sign function $s(\vec{A})\vec{b}$ may be approximated as 
\begin{equation*}
    s(\vec{A})\vec{b} \approx r_q(\vec{A})\vec{b} = \sum_{i=1}^{q} \omega_i \vec{A} (\vec{A}^2 + z_i^2 \vec{I})^{-1} \vec{b}
\end{equation*}
where $z_i$ and $\omega_i$ are appropriately chosen quadrature nodes and weights \cite{hale_higham_trefethen_08}.

We can of course write $r_q(x) = M_q(x) / N_q(x)$, so it is tempting to set $R_q(x) = 1$ and $\vec{H}_q = N_q(\vec{A})$ and then use Lanczos-OR to return the $\vec{H}_q$-norm optimal approximation to $r_q(\vec{A})\vec{b}$ as an approximation for $s(\vec{A})\vec{b}$. 
However, while $r_q(x)$ is convergent to $f(x)$ as $q\to\infty$, $N_q(x) := \prod_{i=1}^q (x^2
+ z_i^2 )$ is not convergent to any fixed function.
In fact $N_q(x)$ will increase in degree and $\vec{H}_q$ will be increasingly poorly conditioned.
This presents a numerical difficulty in computing the Lanczos-OR iterate in this limit.
More importantly, it is not clear that it is meaningful to approximate a function in this way. 
Indeed, it seems reasonable to expect that, for fixed $k$, as $q\to\infty$,  our approximation should be convergent to something. 
However, we cannot guarantee $\lanopt_k(r_q,1)$ is convergent in this limit.

Another option is to compute the \emph{term-wise optimal} approximations to each term in the sum representation of $r_q$ and output
\begin{equation*}
    \sum_{i=1}^{q} \omega_i \vec{Q} ([\widehat{\vec{T}}^2]_{:k,:k} + z_i^2\vec{I})^{-1} \vec{T} \vec{e}_0.
\end{equation*}

Interestingly, this is exactly what would be obtained by using Lanczos-OR to approximate each of the corresponding linear systems in the partial fractions decomposition (which is equivalent to a special case of QMR on such systems).
\begin{lemma}
Suppose $z\in\R$ and define $r(x)= 1/(x^2+z^2)$, $R^\pm(x) = x\pm \ii z$, and $r^\pm(x) = 1/R^\pm(x)$.
Then,
$
    \lanopt_k(r,1)
    = \frac{1}{2\ii z} \left( \lanopt_k(r^-,R^+) - \lanopt_k(r^+,R^-) \right).
$
\end{lemma}

\begin{proof}
We have that
$\lanopt_k(r^\pm,R^\mp) = \vec{Q}([\widehat{\vec{T}}^2+|z|^2\vec{I}]_{:k,:k})^{-1} [\vec{T} \mp \ii z \vec{I}]_{:k,:k} \vec{e}_0$.
So
$
    \lanopt_k(r^-,R^+) - \lanopt_k(r^+,R^-)
    = 2 \ii z  \vec{Q}([\widehat{\vec{T}}^2+|z|^2\vec{I}]_{:k,:k})^{-1} \vec{e}_0
    = 2 \ii z \lanopt_k(r,1).
$
The result follows by rearranging the previous expression.
\end{proof}

Whether it is better to use Lanczos-OR with $r(x)$ and $R(x)=1$ or with $r^\pm(x)$ and $R^\pm(x)$ (i.e. QMR) is somewhat unclear.
Lanczos-OR avoids the need for complex arithmetic, which simplifies implementation slightly. 
However, since QMR has been studied longer, it is likely to have more practical low-memory implementations.

\subsection{An induced approximation}
\label{sec:sign_function}
The approach from the previous section can be taken a step further to obtain from Lanczos-OR what we called an ``induced'' approximation for functions like the sign function. Instead of discretizing an integral representation of the function, we can use it directly. In particular,  
for any $a>0$,
$
    {1}/{\sqrt{a}}
    = \frac{2}{\pi} \int_{0}^{\infty} \frac{1}{a + z^2} \d{z}.
$
Thus, if $s(x) = \operatorname{sign}(x)=x/|x| = x/\sqrt{x^2}$, we have 
\begin{equation*}
    s(\vec{A})\vec{b} = \frac{2}{\pi} \int_{0}^{\infty} \vec{A}(\vec{A}^2 + z^2\vec{I})^{-1} \vec{b} \: \d{z}.
\end{equation*}
The Lanczos-OR approximation to $\vec{A}(\vec{A}^2 + z^2\vec{I})^{-1} \vec{b}$ (with $R(x) = 1$) is $\vec{Q} ([\widehat{\vec{T}}^2]_{:k,:k} + z^2\vec{I})^{-1} \vec{T} \vec{e}_0$, which is optimal over the Krylov subspace in the $(\vec{A}^2 + z^2\vec{I})$-norm. Plugging this approximation into the integral above
 yields the approximation 
\begin{equation*}
    \frac{2}{\pi} \int_{0}^{\infty}  \vec{Q} ([\widehat{\vec{T}}^2]_{:k,:k} + z^2\vec{I})^{-1} \vec{T} \vec{e}_0 \: \d{z}
    = \vec{Q} \left([\widehat{\vec{T}}^2]_{:k,:k} \right)^{-1/2} \vec{T}  \vec{e}_0.
\end{equation*}
Thus, we can define an induced iterate for the matrix sign function as
\begin{equation*}
\textup{\textsf{sign-OR}}_k := \vec{Q} \left([\widehat{\vec{T}}^2]_{:k,:k} \right)^{-1/2} \vec{T}  \vec{e}_0
=
\vec{Q} \left( [\widehat{\vec{T}}]_{:k,:k+1}[\widehat{\vec{T}}]_{:k+1,:k} \right)^{-1/2} \vec{T}  \vec{e}_0.
\end{equation*}
As seen in \cref{sec:experiment_sign_function}, this Lanczos-OR induced iterate, $\textup{\textsf{sign-OR}}_k$, performs very well empirically, and in fact appears to provide close to an optimal approximation from the Krylov subspace for our test problem. It outperforms the standard Lanczos-FA algorithm, exhibiting smoother convergence. However, both methods appear to converge at roughly the same overall rate, and perform remarkably close to optimally. In the following subsection we seek to explain why the iterates behave similarly.

\subsubsection{Relation to Lanczos-FA}
The standard Lanczos-FA method for matrix function approximation is defined as follows:
\begin{definition}
\label{def:lanczos_fa}
The Lanczos-FA iterate is defined as
\begin{equation*}
    \lan_k(f)
    :=  \vec{Q} f(\vec{T}) \vec{e}_0.
 \end{equation*}
\end{definition}
For a rational function $r$ and $\tilde{N}, \tilde{M}$ as in \cref{def:lanczos_or}, we have that $\lan_k(r) = \vec{Q}  \tilde{N}(\vec{T})^{-1}\tilde{M}(\vec{T})\vec{e}_0$. This compares to $([\tilde{N}(\widehat{\vec{T}})]_{:k,:k})^{-1} [\tilde{M}(\widehat{\vec{T}})]_{:k,:k} \vec{e}_0$ for the Lanczos-OR iterate. The two expressions are clearly related since $\tilde{N}(\vec{T})$ and $[\tilde{N}(\widehat{\vec{T}})]_{:k,:k}$ differ only in the bottom rightmost $(q-1)\times (q-1)$ principle submatrix, where $q = \deg(\tilde{N})$.

Using this fact, it can be argued that the Lanczos-OR and Lanczos-FA iterates ``tend to coalesce as convergence takes place'' \cite[Proposition 5.1]{lopez_simoncini_06}.
We show that a similar phenomenon occurs with the induced Lanczos-OR approximation to the sign function and the Lanczos-FA approximation.

\begin{theorem}
\label{thm:sign_approx_compare}
Let $\sigma_{\textup{max}}(\vec{T})$ and $\sigma_{\textup{min}}(\vec{T})$ be the largest and smallest singular values of $\vec{T}$, respectively.
The Lanczos-FA and induced Lanczos-OR approximations to the matrix sign function satisfy
\begin{equation*}
    \| \lan_k(\operatorname{sign}) - \textup{\textsf{sign-OR}}_k \|_2
    \leq
    \tfrac{1}{2}(\beta_{k-1})^2 
    {\sigma_{\textup{max}}(\vec{T})}/{ \sigma_{\textup{min}}(\vec{T})^3}.
\end{equation*}
\end{theorem}

\begin{proof}
We proceed similarly to the proof of \cite[Proposition 5.1]{lopez_simoncini_06}.
Let $r(x) = x/(x^2+z^2)$ so that ${N}(x) = x^2 + z^2$ and ${M}(x) = x$.
Note that ${N}(\vec{T}) = [{N}(\widehat{\vec{T}})]_{:k,:k} - \beta_{k-1}^2 \vec{e}_{k-1}\vec{e}_{k-1}^\T$ so,
\begin{equation*}
    \vec{T} \vec{e}_0
    = {N}(\vec{T}) {N}(\vec{T}) ^{-1} \vec{T} \vec{e}_0
    = ([{N}(\widehat{\vec{T}})]_{:k,:k} - \beta_{k-1}^2 \vec{e}_{k-1}\vec{e}_{k-1}^\T){N}(\vec{T})^{-1} \vec{T} \vec{e}_0.
\end{equation*}
Thus, left multiplying by $\vec{Q}([{N}(\widehat{\vec{T}})]_{:k,:k})^{-1}$ and rearranging terms we find that
\begin{align*}
    \lan_k(r) - \lanopt_k(r,1)
    = \beta_{k-1}^2 \vec{Q}([{N}(\widehat{\vec{T}})]_{:k,:k})^{-1} \vec{e}_{k-1}\vec{e}_{k-1}^\T {N}(\vec{T})^{-1} \vec{T} \vec{e}_0.
\end{align*}
Now, suppose that $f(x) = \operatorname{sign}(x)$ and set $ \textsf{error}_k := \lan_k(\operatorname{sign}) - \textup{\textsf{sign-OR}}_k$.
Then, since the Lanczos-FA approximation can also be induced by an integral over $z\in[0,\infty)$, we have that,
\begin{equation*}
    \textsf{error}_k
    = \beta_{k-1}^2 
    \frac{2 }{\pi} \int_0^\infty
    \vec{Q} ([{N}(\widehat{\vec{T}})]_{:k,:k})^{-1} \vec{e}_{k-1}\vec{e}_{k-1}^\T (\vec{T}^2 + z^2 \vec{I})^{-1} \vec{T} \vec{e}_0 
    \d{z}.
\end{equation*}
Note that $[\widehat{\vec{T}}^2]_{:k,:k} - \vec{T}^2 = \beta_{k-1}^2 \vec{e}_{k-1} \vec{e}_{k-1}^\T$ is positive semidefinite.
Therefore, using that $\sigma_{\textup{min}}([\widehat{\vec{T}}^2]_{:k,:k}) \geq \sigma_{\textup{min}}(\vec{T}^2) = \sigma_{\textup{min}}(\vec{T})^2$, and $[N(\widehat{\vec{T}})]_{:k,:k} = [\widehat{\vec{T}}^2]_{:k,:k} + z^2 \vec{I}$,
\begin{align*}
    \| \textsf{error}_k \|_2
    &
    = \beta_{k-1}^2 
     \left\| \vec{Q} \left( \frac{2 }{\pi} \int_0^\infty
    ([\widehat{\vec{T}}^2]_{:k,:k} + z^2 \vec{I})^{-1} \vec{e}_{k-1}\vec{e}_{k-1}^\T (\vec{T}^2 + z^2 \vec{I})^{-1} \d{z} \right)  \vec{T} \vec{e}_0 
     \right\|_2
    \\&
    \leq \beta_{k-1}^2 
     \left( \frac{2 }{\pi} \int_0^\infty
    \| ([\widehat{\vec{T}}^2]_{:k,:k} + z^2 \vec{I})^{-1} \|_2 \| (\vec{T}^2 + z^2 \vec{I})^{-1} \|_2 \d{z} \right)  \| \vec{T} \|_2
    \\& \leq
    \beta_{k-1}^2 
     \left( \frac{2 }{\pi} \int_0^\infty
    | (\sigma_{\textup{min}}(\vec{T})^2 + z^2 )^{-1} |~ | (\sigma_{\textup{min}}(\vec{T})^2 + z^2 )^{-1} | \d{z} \right)  \sigma_{\textup{max}}(\vec{T})
    \\&= \beta_{k-1}^2
    \frac{\sigma_{\textup{max}}(\vec{T})}{2 \sigma_{\textup{min}}(\vec{T})^3}.
\end{align*}
\end{proof}

Since $|\beta_{k-1}|$ tends to decrease as the Lanczos method converges, this seemingly implies that the induced Lanczos-OR iterate and the Lanczos-FA iterate tend to converge in this limit.
However, recall that $\vec{T} = [\widehat{\vec{T}}]_{:k,:k}$ changes at each iteration $k$. 
Thus, there is the difficulty that $\vec{T}$ may have an eigenvalue near zero, in which case the preceding bound could be useless. 
However, it is known that $\vec{T}$ cannot have eigenvalues near zero in two successive iterations, assuming that the eigenvalues of $\vec{A}$ are not too close to zero.
Specifically,  \cite[Equation 3.10]{greenbaum_druskin_knizhnerman_99} asserts that
\begin{equation*}
\max\{ \sigma_{\textup{min}}([\vec{T}]_{:k-1}),\sigma_{\textup{min}}([\vec{T}]_{:k,:k})\} > \frac{\sigma_{\textup{min}}(\vec{A})^2}{(2+\sqrt{3})\|\vec{A}\|_2}.
\end{equation*}
Since $\beta_{k-1}$ has little to do with the minimum magnitude eigenvalue of $\vec{T}$ (recall that the Lanczos recurrence is shift invariant), we expect that the induced Lanczos-OR iterate and the Lanczos-FA iterate will become close as the Lanczos algorithm converges, at least at every other iteration.
This implies that existing spectrum dependent bounds for Lanczos-FA for the sign function \cite{chen_greenbaum_musco_musco_22} can be carried over to Lanczos-OR.
More interestingly, it means that understanding the induced approximation may provide a way of understanding Lanczos-FA for the matrix sign function.
Since Lanczos-FA often exhibits oscillatory behavior, bounds for the induced Lanczos-OR based approximation may be easier to obtain.

\section{Implementing Lanczos-OR using low memory}
\label{sec:low_mem}

We now describe a low-memory implementation of Lanczos-OR which is similar in spirit to CG, MINRES, and QMR. 
It is inspired by the LDL based version of CG described in \cite{simonova_tichy_21} and is closely related to the DIOM method in \cite[Section 6.4]{saad_03}.
A full NumPy implementation, including the code required to reproduce all of our experiments, is available online.

For convenience, we will denote $\vec{M} := [\tilde{M}(\widehat{\vec{T}})]_{:k,:k}$ and $\vec{N} := [\tilde{N}(\widehat{\vec{T}})]_{:k,:k}$. Thus, the Lanczos-OR output is given by $\vec{Q} \vec{N}^{-1} \vec{M} \vec{e}_0$.
At a high level, our approach is to:
\begin{itemize}
    \item Take one iteration of Lanczos to generate one more column of $\widehat{\vec{Q}}$ and $\widehat{\vec{T}}$
    \item Compute one more column of $\vec{M}$ and $\vec{N}$
    \item Compute one more factor of $\vec{L}^{-1} = \vec{L}_{k-1}\cdots \vec{L}_1 \vec{L}_0$ and one more entry of $\vec{D}$ where $\vec{L}$ and $\vec{D}$ are defined by the LDL factorization $\vec{N} =\vec{L}\vec{D}\vec{L}^\T$
    \item Compute one more term of the sum:
    \begin{equation*}
        \vec{Q}\vec{N}^{-1}\vec{M}\vec{e}_0
        = \vec{Q} \vec{L}^{-1} \vec{D}^{-1} \vec{L}^{-\T} \vec{M} \vec{e}_0 
        = \sum_{i=0}^{k-1} \frac{[\vec{L}^{-\T}\vec{M}\vec{e}_0]_i}{[\vec{D}]_{i,i}} [\vec{Q}\vec{L}^{-1}]_{:,i}
    \end{equation*}
\end{itemize}

There are two critical observations which must be made in order to see that this gives a memory-efficient implementation.
The first is that, since $\widehat{\vec{T}}$ is tridiagonal, $\vec{M}$, $\vec{N}$, and therefore $\vec{L}$ are all of half-bandwidth $q := \max(\deg(\tilde{M}),\deg(\tilde{N}))$.
This means that it is possible to compute the entries of $\vec{D}$ and the factors of $\vec{L}^{-1} = \vec{L}_{k-1}\cdots \vec{L}_1 \vec{L}_0$ one by one as we get the entries of $\widehat{\vec{T}}$.
The second is that because $\vec{L}$ is of bandwidth $q$, we can compute $[\vec{Q}\vec{L}^{-1}]_{:,i}$ without saving all of $\vec{Q}$. 
More specifically, $[\vec{L}^{-1} \vec{M} \vec{e}_0]_i$ and $[\vec{Q} \vec{L}^{-1}]_{:,i}$ can respectively be computed from $\vec{L}_{j-1}\cdots \vec{L}_1 \vec{L}_0 \vec{M} \vec{e}_0$ and $\vec{Q} \vec{L}_0^\T \vec{L}_1^\T \cdots \vec{L}_{k-1}^\T$ and can therefore be maintained iteratively as the factors of $\vec{L}^{-1}$ are computed. 
Moreover, because of the banded structure of the factors $\vec{L}_i$, we need only maintain a sliding window of the columns of $\vec{Q}\vec{L}^{-1}$ which will allow us to access the relevant columns when we need them and discard them afterwards.

The cost of such an implementation of Lanczos-OR is $O((k+q)(T_{\textup{mv}}+n))$, where $T_{\textup{mv}}$ is the cost of a matrix-vector product.
On the other hand, Lanczos-FA, implemented by a similar LDL factorization of $\vec{T}$ would require $O(k(T_{\textup{mv}}+n))$.
Since $q$ is a constant typcially far smaller than $k$, this is unlikely to be of major concern.
For instance, in many of our numerical experiments $q=1$ while $k>100$.

We now provide the details of the implementation.
For clarity, we only describe how to compute $\vec{M}$ and $\vec{N}$ in the case that $\tilde{M}(x)$ and $\tilde{N}(x)$ are degree at most two.
The rest of the subroutines are fully described for any degree.
The syntax we use follows Python and other object oriented languages closely.

\subsection{Computing LDL factorization}

For the time being, assume that we can sequentially access the rows of $\vec{M}$ and $\vec{N}$.
Our first step is to compute an LDL factorization of $\vec{N}$, which can be done using a symmetrized version of Gaussian elimination and is guaranteed to exist if \( \vec{N} \) is positive definite \cite{higham_02}.
Specifically, Gaussian elimination can be viewed as transforming the starting matrix \( \vec{N}_0 = \vec{N} \) to a diagonal matrix \( \vec{N}_{k-1} = \vec{D} \) via a sequence of row and column operations
$
    \vec{N}_{i+1} = \vec{L}_i \vec{N}_i \vec{L}_i^\T,
$
where
\begin{equation*}
    \vec{L}_i := 
    \vec{I}_k + \vec{l}_i\vec{e}_i^\T
    ,\qquad
    \vec{l}_i := 
        -\big[ 0 , \cdots , 0 , {[\vec{N}_{i}]_{i+1,i}}/{[\vec{N}_i]_{i,i}}  , \cdots , {[\vec{N}_i]_{k-1,i}}/{[\vec{N}_i]_{i,i}}  \big]^\T.
\end{equation*}
Note that the entries of \( \vec{L}_i \) are chosen to introduce zeros to the \( i \)-th row and column of \( \vec{N}_i \) such that \( [\vec{N}_{i+1}]_{:i+1,:i+1} \) is diagonal.
Therefore, if the algorithm terminates successfully, we will have obtained a factorization 
\begin{equation*}
    \vec{D} = (\vec{L}_{k-1} \cdots \vec{L}_1 \vec{L}_0) \vec{N} (\vec{L}_0^\T \vec{L}_1^\T \cdots \vec{L}_{n-1}^\T)
\end{equation*}
where \( \vec{D} \) is diagonal and each \( \vec{L}_i \) is unit lower triangular.
To obtain the factorization \( \vec{N} = \vec{L} \vec{D} \vec{L}^{\T} \), simply define  \( \vec{L} := (\vec{L}_{k-1} \cdots \vec{L}_1 \vec{L}_0)^{-1} \) and note that 
\(
    \vec{L} = \vec{I}_k - \sum_{i=0}^{k-1} \vec{l}_i \vec{e}_i^\T.
\)
Observe that $\vec{l}_{k-1}$ is the zeros vector and is only included in sums for ease of indexing later on. 
For further details on LDL factorizations, we refer readers to \cite{higham_02}.
To implement a LDL factorization, observe that the procedure above defines a recurrence
\begin{equation*}
    [\vec{D}]_{j,j} = [\vec{N}]_{j,j} - \sum_{\ell=0}^{j-1} [\vec{L}_{j,\ell}]^2 [\vec{D}]_{\ell,\ell}
    ,\quad
    [\vec{L}]_{i,j} = \frac{1}{[\vec{D}]_{j,j}} \bigg( [\vec{N}]_{i,j} - \sum_{\ell=0}^{j-1}[\vec{L}]_{j,\ell} [\vec{L}]_{i,\ell} [\vec{D}]_{\ell,\ell} \bigg).
\end{equation*}

The fact that $\vec{L}$ has the same half bandwidth as $\vec{N}$ allows an efficient LDL implementation where terms which are known to be zero are not computed and only the important diagonals of $\vec{L}$ are stored.
This implementation is fed a stream of the columns of $\vec{N}$ in order, as shown in \cref{fig:streaming_LDL}.
Here the diagonal of $\vec{D}$ is stored as $\ms{d}$ and the $(j+1)$-st diagonal of $\vec{L}$ is stored as $[\ms{L}]_{j,:}$.
Thus, $\vec{L}_{i,j} = [\ms{L}]_{i-j-1,j}$ as long as $i-j \in \{0,1,\ldots, q\}$.
Note that this implementation is equivalent, even in finite precision arithmetic, to the standard implementation based on the above recurrences.

\begin{algorithm}[htb]
\caption{Streaming LDL}\label{alg:streaming_LDL}
\fontsize{10}{10}\selectfont
\begin{algorithmic}[1]
\Class{streaming-LDL}{$q,k$}
\State{stream: $[\vec{N}]_{0,0:q+1}, [\vec{N}]_{1,1:q+2}, \ldots, [\vec{N}]_{k-1,k-1:q+k-1}$}
\State $\ms{L} = \textsc{zeros}(q,k)$
\State $\ms{d} = \textsc{zeros}(k)$
\State $\ms{j}\gets0$
\Procedure{read-stream}{$\vec{n}$}
\State $[\ms{d}]_{\ms{j}} \gets [\vec{n}]_0 - \sum_{\ell=\max(0,\ms{j}-q)}^{\ms{j}-1} [\ms{L}]_{\ms{j}-\ell-1,\ell}^2 [\ms{d}]_{\ell}$
\For{$i=\ms{j+1}, \ldots,  \min(\ms{j}-q,k-1)$}
\State $[\ms{L}]_{i-\ms{j}-1,\ms{j}} \gets (1/[\ms{d}]_{\ms{j}})([\vec{n}]_{i-\ms{j}} - \sum_{\ell=\max(0,i-q)}^{i-1} [\ms{L}]_{i-\ell-1,\ell} [\ms{L}]_{\ms{j}-\ell-1,\ell} [\ms{d}]_{\ell} )$
\EndFor
\State $\ms{j}\gets \ms{j}+1$
\EndProcedure
\EndClass
\end{algorithmic}
\end{algorithm}

\subsection{Inverting the LDL factorization}

Once we have computed a factorization \( \vec{N} = \vec{L}\vec{D}\vec{L}^\T \), we can easily evaluate $\vec{Q}\vec{L}^{-1}\vec{D}^{-1} \vec{L}^{-\T}\vec{M}\vec{e}_1$ using the fact that \( \vec{L}^{-1} = \vec{L}_{k-1} \cdots \vec{L}_1 \vec{L}_0 \).
Moreover, because the \( \vec{L}_j \) can be computed one at a time, there is hope that we can derive a memory efficient implementation.

Towards this end, define \( \vec{y}_j := \vec{L}_{j-1} \cdots \vec{L}_1 \vec{L}_0 \vec{M}\vec{e}_0 \) and \( \vec{X}_j := \vec{Q} \vec{L}_0^\T \vec{L}_1^\T \cdots \vec{L}_{j-1}^\T  \).
Then, setting \( \vec{y}_{0} = \vec{M} \vec{e}_1 \) we have that
\begin{equation*}
    \vec{y}_{j+1} 
    = \vec{L}_{j} \vec{y}_{j}
    = (\vec{I} + \vec{l}_j \vec{e}_j^\T ) \vec{y}_{j}
    = \vec{y}_{j} +  (\vec{e}_{j}^\T \vec{y}_{j}) \vec{l}_j.
\end{equation*}
Similarly, setting \( \vec{X}_{0} = \vec{Q} \) we have that
\begin{equation*}
    \vec{X}_{j+1} 
    = \vec{X}_{j} \vec{L}_j^\T
    = \vec{X}_{j} ( \vec{I} + \vec{e}_j \vec{l}_j^\T)
    = \vec{X}_{j} + \vec{X}_{j} \vec{e}_j \vec{l}_j^\T.
\end{equation*}
Then \( \vec{Q}\vec{L}^{-1}\vec{D}^{-1} \vec{L}^{-\T}\vec{M}\vec{e}_1 = \vec{X}_{k} \vec{D}^{-1} \vec{y}_{k} \) can be computed accessing \( \vec{L} \), and therefore \( \vec{N} \), column by column.

\begin{algorithm}[htb]
\caption{Streaming banded product}\label{alg:streaming_banded_prod}
\fontsize{10}{10}\selectfont
\begin{algorithmic}[1]
\Class{streaming-banded-prod}{$n,k,q$}
\State{stream:$ $}
\State $\ms{X\_} \gets \textsc{zeros}(n,q+1)$
\State $\ms{y\_} \gets \textsc{zeros}(q+1)$
\State $\ms{out} \gets \textsc{zeros}(n)$
\State $\ms{j}\gets-1$
\Procedure{read-stream}{$\vec{v},\vec{l},d,\vec{y}_0$}
\If{$\ms{j}=-1$}
\State $[\ms{X\_}]_{:,:q} = \vec{v}$
\Else
\If{$\ms{i}=-1$}
\State $\ms{y\_} \gets \vec{y}_0$
\EndIf
\State $\ms{out} \gets \ms{out} + ([\ms{y\_}]_{0}/d)[\ms{X\_}]_{:,0}$
\State $[\ms{y\_}]_{:q} \gets [\ms{y\_}]_{1:} - [\ms{y\_}]_{0} \vec{l}$
\State $[\ms{y\_}]_{-1} \gets 0$
\State $[\ms{X\_}]_{:,-1} \gets \vec{v}$
\State $[\ms{X\_}]_{:,:q} \gets [\ms{X\_}]_{:,1:} + [\ms{X\_}]_{:,0} \vec{l}^\T$
\EndIf
\State $\ms{j}\gets \ms{j}+1$
\EndProcedure
\EndClass
\end{algorithmic}
\end{algorithm}

\subsubsection{Streaming version}
\label{sec:streaming_banded_solve}

Recall that $[\vec{l}_i]_{:,\ell}$ is zero if $\ell \leq i$ or $\ell > i +q$.
Since $[\vec{l}_i]_{:i}$ is zero, we have
\begin{align*}
    [\vec{y}_j]_j
    &= [\vec{y}_{j} +  (\vec{e}_{j}^\T \vec{y}_{j}) \vec{l}_j]_j
    = [\vec{y}_{j+1}]_j
    =\cdots
    = [\vec{y}_{k}]_j
\\
    [\vec{X}_{j}]_{:,j}
    &= [\vec{X}_{j} + \vec{X}_{j} \vec{e}_j \vec{l}_j^\T]_{:,j}
    = [\vec{X}_{j+1}]_{:,j}
    = \cdots = [\vec{X}_{k}]_{:,j}.
\end{align*}
We therefore have that
\begin{align*}
    \vec{X}_{k} \vec{D}^{-1} \vec{y}_{k}
    = \sum_{j=0}^{k-1} \frac{[\vec{y}_{k}]_j}{[\vec{D}]_{j,j}} [\vec{X}_{k}]_{:,j}
    = \sum_{j=0}^{k-1} \frac{[\vec{y}_{j}]_j}{[\vec{D}]_{j,j}} [\vec{X}_{j}]_{:,j}.
\end{align*}
Similarly, since $[\vec{l}_i]_{i+q+1:}$ is zero,
\begin{align*}
    [\vec{y}_j]_{j+q:}
    &= [\vec{y}_{j-1} + (\vec{e}_{j-1}^\T\vec{y}_{j-1}) \vec{l}_{j-1}]_{j+q:}
    = [\vec{y}_{j-1}]_{j+q:}
    = \cdots
    = [\vec{y}_0]_{j+q:}
    \\
    [\vec{X}_{j}]_{:,j+q:}
    &= [\vec{X}_{j-1}+\vec{X}_{j-1}\vec{e}_{j-1}\vec{l}_{j-1}^\T]_{:,j+q:}
    = [\vec{X}_{j-1}]_{:,j+q:}
    =\cdots
    = [\vec{X}_{0}]_{:,j+q:}.
\end{align*}
By definition, $\vec{y}_0 = \vec{M}\vec{e}_0$ and $\vec{X}_0 = \vec{Q}$.
Thus, we see that it is not necessary to know the later columns of $\vec{X}_j$ immediately.

\begin{algorithm}[htb]
\caption{Streaming banded inverse}\label{alg:streaming_banded_inv}
\fontsize{10}{10}\selectfont
\begin{algorithmic}[1]
\Class{streaming-banded-inv}{$n,k,q$}
\State{stream:$ $}
\State $\ms{LDL} \gets \textsc{streaming-LDL}(k,q)$
\State $\ms{Q0} \gets \textsc{zeros}(n,q)$
\State $\ms{j}\gets0$
\Procedure{read-stream}{$\vec{q},\vec{n},\vec{y}_0$}
\If{$\ms{j}<q$}
\State $[\ms{Q0}]_{:,j} \gets \vec{v}$
\If{$\ms{j}=q-1$}
    \State $\ms{b-prod} \gets \textsc{streaming-banded-prod}(n,k,q)$
    \State $\ms{b-prod.}\textsc{read-stream}(\ms{V0},\ms{none},\ms{none},\ms{none})$
\EndIf
\Else
\State $\ms{LDL.}\textsc{read-stream}(\vec{n})$
\State $\ms{b-inv.}\textsc{read-stream}(\vec{q},-[\ms{LDL.}\ms{L}]_{:,j-q},[\ms{LDL.}\ms{d}]_{j-q},\vec{y}_0)$
\EndIf
\State $\ms{j}\gets \ms{j}+1$
\EndProcedure
\EndClass
\end{algorithmic}
\end{algorithm}

We can define a streaming algorithm by maintaining only the relevant portions of the $\vec{X}_i$ and $\vec{y}_i$.
Towards this end, define the length \( q+1 \) vector \( \bar{\vec{y}}_j := [\vec{y}_j]_{j:j+q+1} \) the \( n\times (q+1) \) matrix \( \bar{\vec{X}}_j = [\vec{X}_j]_{:,j:j+q+1} \).
Using the above observations, we see that these quantities can be maintained by the recurrences
\begin{align*}
    \bar{\vec{y}}_j &=  
    \begin{bmatrix}
        [\bar{\vec{y}}_{j-1}]_{1:}\\0
    \end{bmatrix}
    +  [\bar{\vec{y}}_{j-1}]_0 [\vec{l}_j]_{j+1:j+q+1}
    \\
    [\bar{\vec{X}}_j]_{:,:q}
    &= [\bar{\vec{X}}_{j-1}]_{:,1:} 
        + ( [\bar{\vec{X}}_{j-1}]_{:,1} ) ([\vec{l}_j]_{j+1:j+q+1})^\T,
    \qquad
    [\bar{\vec{X}}_j]_{:,q} = [\vec{Q}]_{j+q}.
\end{align*}
Note then that,
\begin{equation*}
    \vec{X}_{k-1} \vec{D}^{-1} \vec{y}_{k-1}
    = \sum_{j=0}^{k-1}  \frac{[\bar{\vec{y}}_j]_{1}}{[\vec{D}]_{j+1,j+1}} [\bar{\vec{X}}_j]_{:,1}.
\end{equation*}
This results in \cref{alg:streaming_banded_prod,alg:streaming_banded_inv}, whose streaming data access patterns are outlined in \cref{fig:streaming_banded_prod}.

\begin{algorithm}[htb]
\caption{Streaming tridiagonal square}\label{alg:streaming_tridiag_square}
\fontsize{10}{10}\selectfont
\begin{algorithmic}[1]
\Class{streaming-tridiagonal-square}{$k$}
\State{stream: $(\alpha_0, \beta_0), \ldots, (\alpha_{k-1},\beta_{k-1})$}
\State $\ms{T} \gets \textsc{zeros}(2,k)$
\State $\ms{Tp2} \gets \textsc{zeros}(3,k)$
\State $\ms{j}\gets0$
\Procedure{read-stream}{$\alpha,\beta$}
\State $[\ms{T}]_{0,\ms{j}} = \alpha$
\State $[\ms{T}]_{1,\ms{j}} = \beta$
\If{$\ms{i}=0$}
\State $[\ms{Tp2}]_{0,\ms{j}} \gets [\ms{T}]_{0,\ms{j}}^2+[\ms{T}]_{1,\ms{j}}^2$
\Else
\State $[\ms{Tp2}]_{0,\ms{j}} \gets [\ms{T}]_{0,\ms{j}}^2+[\ms{T}]_{1,\ms{j}}^2+[\ms{T}]_{1,\ms{j}-1}^2$
\State $[\ms{Tp2}]_{1,\ms{j}} \gets ([\ms{T}]_{0,\ms{j}}+[\ms{T}]_{0,\ms{j}-1})[\ms{T}]_{1,\ms{j}-1}$
\State $[\ms{Tp2}]_{2,\ms{j}} \gets [\ms{T}]_{1,\ms{j}}[\ms{T}]_{1,\ms{j}-1}$
\EndIf
\State $\ms{j}\gets \ms{j}+1$
\EndProcedure
\EndClass
\end{algorithmic}
\end{algorithm}

\begin{algorithm}[htb]
\caption{Get polynomial of tridiagonal matrix}\label{alg:get_poly}
\fontsize{10}{10}\selectfont
\begin{algorithmic}[1]
\Procedure{get-poly}{$P,\ms{STp2},k,j$}
\State $a,b,c = P(0), P'(0), P''(0)$
\State $\vec{p} \gets \textsc{zeros}(3)$
\State $[\vec{p}]_{:3} \gets a [\ms{STp2.}\ms{Tp2}]_{:,j}$
\State $[\vec{p}]_{:2} \gets b [\ms{STp2.}\ms{T}]_{:,j}$
\State $[\vec{p}]_{:1} \gets c $
\EndProcedure
\end{algorithmic}
\end{algorithm}

\subsection{Computing polynomials in $\vec{T}$}

The last major remaining piece is to construct $\vec{M} = \tilde{M}(\vec{T})$ and $\vec{N} = [\tilde{N}(\widehat{\vec{T}})]_{:k,:k}$.
Recall that we have assumed $\tilde{M}$ and $\tilde{N}$ are of degree at most two for convenience.
In iteration $\ell$ of Lanczos, we obtain $\alpha_{\ell}$ and $\beta_{\ell}$.
Observe that $\vec{T}^2$ is symmetric and that, defining $\beta_{-1} = \beta_{k} = 0$, the lower triangle is given by
$[\vec{T}^2]_{i,j} = \beta_{j-1}^2 + \alpha_j^2 + \beta_j^2$ if $j=i$, $[\vec{T}^2]_{i,j} = (\alpha_j + \alpha_{j+1}) \beta_i$ if $j=i-1$, $[\vec{T}^2]_{i,j} = \beta_j \beta_{j+1}$ if $j=i-2$, and $0$ elsewhere.

\begin{figure}[htb]
    \centering
    \subfloat[\textit{Pattern for $\vec{N}$ in \cref{alg:streaming_LDL}.}]{%
    \includegraphics[height=2.70cm]{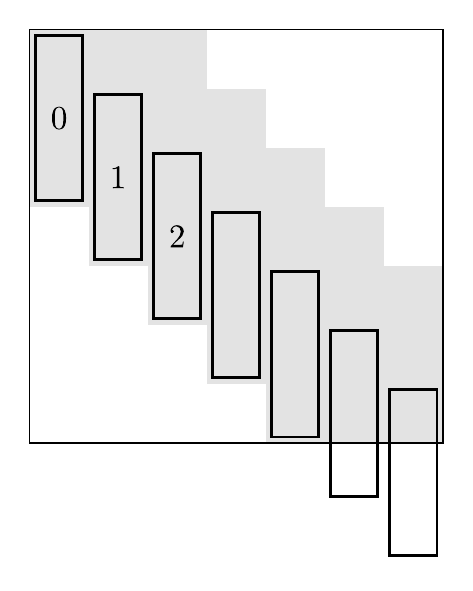}
    \label{fig:streaming_LDL}
    }
    \hfill
    \subfloat[\textit{Pattern for $\vec{Q}$, $\vec{L}$, $\vec{d}$, $\vec{M}$ in \cref{alg:streaming_banded_prod}.}]{%
    \includegraphics[height=2.75cm]{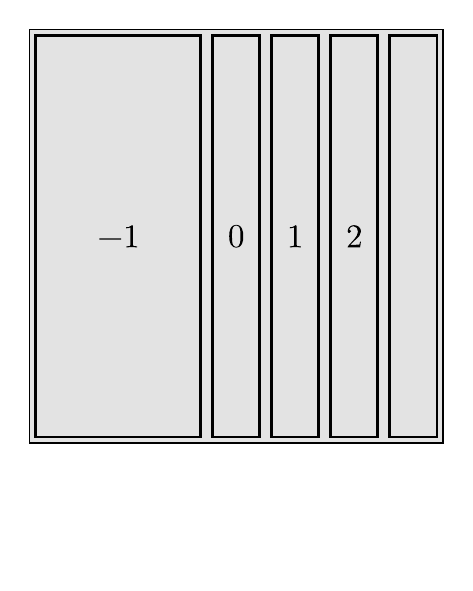}
    \hfill   
    \includegraphics[height=2.75cm]{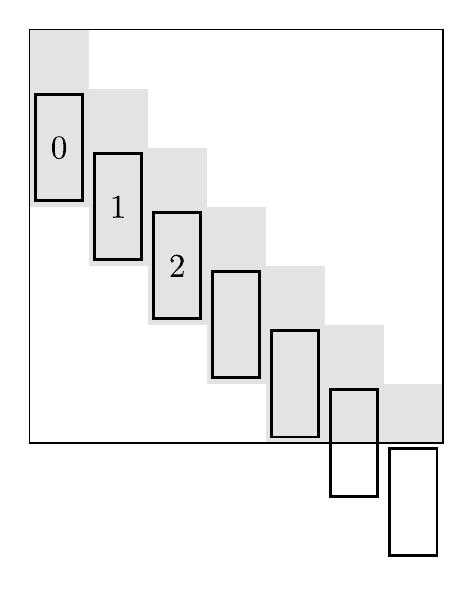}
    \hfill    
    \includegraphics[height=2.75cm]{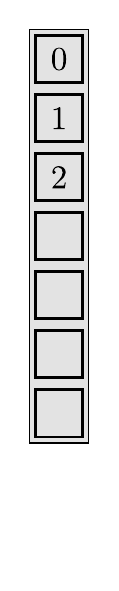}
    \hfill
    \includegraphics[height=2.75cm]{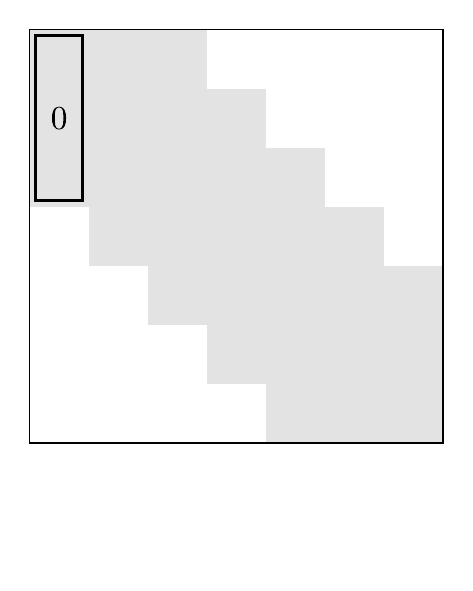}
    \label{fig:streaming_banded_prod,alg:streaming_banded_inv}
    }
    \hfill
    \subfloat[\textit{Pattern for $\vec{T}$ in \cref{alg:streaming_tridiag_square}.}]{%
    \includegraphics[height=2.75cm]{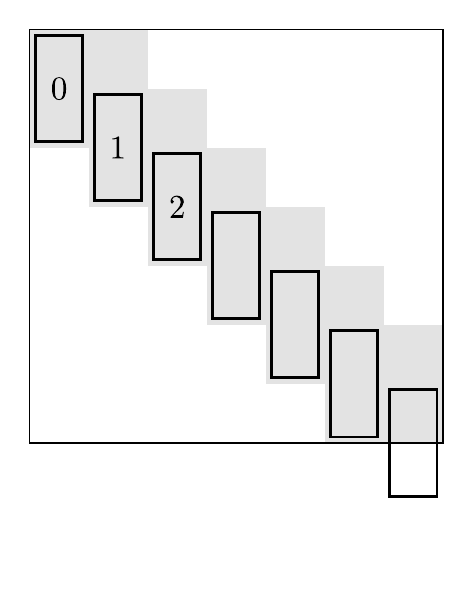}
    \label{fig:streaming_tridiag_square}
    }
    \caption{Access patterns for inputs to streaming functions used in low-memory implementations of Lanczos-OR and Lanczos-FA.
    Indices indicate what information should be streamed into the algorithm at the given iteration.}
   \label{fig:streaming_banded_prod}
\end{figure}

We can use this to implement the streaming algorithm, \cref{alg:streaming_tridiag_square}, for computing the entries of $\vec{T}^2$.
Rather than being fed the entire tridiagonal matrix $\vec{T}$, \cref{alg:streaming_tridiag_square} is fed a stream of the columns of $\vec{T}$ in order, as shown in \cref{fig:streaming_tridiag_square}.
The algorithm respectively stores the $j$-th diagonals of $\vec{T}$ and $\vec{T}^2$  as $[\ms{T}]_{j,:}$ and $[\ms{Tp2}]_{j,:}$.
Then, since we maintain the columns of $\vec{T}^2$ with \cref{alg:streaming_tridiag_square}, we can easily compute $\vec{M}$ and $\vec{N}$ using \cref{alg:get_poly}.

\subsection{Putting it all together}

With this algorithm in place, putting everything together is straightforward, and the full implementation is shown in \cref{alg:rational_lanczos}.
This can be incorporated into any Lanczos implementation and used to compute the Lanczos-OR iterates.
For concreteness, we show this with a standard implementation of Lanczos.
We call the resulting implementation Lanczos-OR-lm.

We can easily obtain an implementation of Lanczos-FA, which we call Lanczos-FA-lm, by replacing $\beta_{k-1}$ with $0$ in the final iteration of the loop.

\begin{algorithm}[htb]
\caption{Streaming banded rational inverse}\label{alg:rational_lanczos}
\fontsize{10}{10}\selectfont
\begin{algorithmic}[1]
\Class{banded-rational}{$n,k,\tilde{M},\tilde{N}$}
\State $\ms{b-inv} \gets \textsc{banded-inv}(n,k,2)$
\State $\ms{STp2} \gets \textsc{streaming-tridiagonal-square}(k)$
\State $\ms{j}\gets0$
\Procedure{read-stream}{$\vec{q},\alpha,\beta$}
\If{$\ms{j}<k$}
    \State $\ms{STp2.}\textsc{read-stream}(\alpha,\beta)$
    \State $\ms{b-inv.}\textsc{read-stream}($
    \State \hspace{1em} $\vec{q},$
    \State \hspace{1em} $\textsc{getpoly}(\tilde{N},\ms{STp2},k,\ms{j}-1)~\textbf{if}~\ms{j}\geq 2~\textbf{else}~\ms{none},$
    \State \hspace{1em} $\textsc{getpoly}(\tilde{M},\ms{STp2},k,\ms{j}-1)~\textbf{if}~\ms{j}= 2~\textbf{else}~\ms{none},$
    \State$)$
\EndIf
\State $\ms{LDL.}\textsc{read-stream}(\vec{n})$
\State $\ms{j}\gets \ms{j}+1$
\EndProcedure
\Procedure{finish-up}{$ $}
\For{$i=k:k+2$}
    \State $\ms{b-inv.}\textsc{read-stream}( \ms{none}, \textsc{getpoly}(\tilde{N},\ms{STp2},k,\ms{j}-1), \ms{none} )$
\EndFor
\EndProcedure
\Procedure{get-output}{$ $}
\State \Return \ms{b-inv.}\ms{b-prod.}\ms{out}
\EndProcedure
\EndClass 
\end{algorithmic}
\end{algorithm}

\begin{algorithm}[htb]
\caption{Lanczos-OR-lm}
 \label{alg:main}
\fontsize{10}{10}\selectfont
\begin{algorithmic}[1]
\Procedure{Lanczos-OR-lm}{$\vec{A},\vec{b},k,M,N,R$}
\State \( \vec{q}_{-1} = \vec{0} \),
\( \beta_{-1} = 0 \),
\( \vec{q}_0  = \vec{b} / \| \vec{b} \| \)
\State Set $\tilde{M}(x) = M(x)R(x)$ and $\tilde{N}(x) = N(x)R(x)$
\State $\ms{lan-lm} \gets \textsc{banded-rational}(n,k,\tilde{M},\tilde{N})$
\For {\( j=0, \ldots, k-1 \)}
    \State \( \tilde{\vec{q}}_{j+1} = \vec{A} \vec{q}_{j} - \beta_{j-1} \vec{q}_{j-1} \)
    \State \( \alpha_j = \langle \tilde{\vec{q}}_{j+1}, \vec{q}_j \rangle \)
    \State \( \tilde{\vec{q}}_{j+1} = \tilde{\vec{q}}_{j+1} - \alpha_j \vec{q}_j \)
    \State \( \beta_{j} = \| \tilde{\vec{q}}_{j+1} \| \)
    \State \( \vec{q}_{j+1} = \tilde{\vec{q}}_{j+1} / \beta_{j} \)
    \State $\ms{lan-lm.}\textsc{read-stream}(\vec{q}_{j},\alpha_{j},\beta_{j})$
\EndFor
\State $\ms{lan-lm.}\textsc{finish-up}()$
\EndProcedure 
\end{algorithmic}
\end{algorithm}

\subsection{Some comments on implementation}

Our main goal is to describe how to implement Lanczos-FA and Lanczos-OR in a way that requires $k$ matrix-vector products and $O(n)$ storage, when $M$ and $N$ are each at most degree two. As mentioned, the approach can be extended to any constant degree. To obtain possibly improved practical performance, it is possible to slightly optimize the storage requirements of our implementation. For example, the implementation described above saves $\vec{T}$, $\vec{T}^2$, $\vec{L}$, and $\vec{d}$, but only accesses a sliding window of these quantities.
We have chosen to save them for convenience since they require only $O(k)$ storage.
However, storing only the relevant information from these quantities would result in an implementation with storage costs independent of the number of iterations $k$.
In this vein, a practical implementation would likely determine $k$ adaptively by monitoring the residual or other measures of the error.

Improvements to the number of vectors of length $n$ may be possible as well, although we expect these would be limited to constant factors. For example, storage could possibly be reduced by incorporating the Lanczos iteration more explicitly with the inversion of the LDL factorization, much like the classical Hestenes and Stiefel implementation of CG \cite{hestenes_stiefel_52}.

\subsection{Lanczos-FA-lm and Lanczos-OR-lm in finite precision arithmetic}
\label{sec:stability}

As with other short-recurrence based Krylov subspace methods, the behavior of Lanczos-FA-lm and Lanzos-OR-lm in finite precision arithmetic may be different than in exact arithmetic.
Fortunately, quite a bit is known about the standard implementation of Lanczos \cite{paige_71,paige_76,paige_80,greenbaum_89,musco_musco_sidford_18}, and we have stated Lanczos-FA-lm and Lanczos-OR-lm in terms of this implementation. 
Knowledge about the standard implementation of Lanczos carries over to Lanczos-FA. 
For instance, assuming $\vec{Q} f(\vec{T}) \vec{e}_0$ is computed accurately from the output of the standard Lanczos algorithm, many error bounds for Lanczos-FA are still applicable \cite{musco_musco_sidford_18,chen_greenbaum_musco_musco_22}.
It is more or less clear that Lanczos-FA-lm and Lanczos-OR-lm will accurately compute the expressions $\vec{Q} N(\vec{T})^{-1} M(\vec{T}) \vec{e}_0$ and $\vec{Q}  ([\tilde{N}(\widehat{\vec{T}})]_{:k,:k})^{-1} [\tilde{M}(\widehat{\vec{T}})]_{:k,:k} \vec{e}_0$ provided that $N(\vec{T})$, $[\tilde{N}(\widehat{\vec{T}})]_{:k,:k}$ are reasonably well conditioned.
Indeed, in practice solving linear systems by symmetric Gaussian elimination is accurate; see for instance \cite[Chapter 10]{higham_02}.
Thus, such bounds and techniques can be applied to Lanczos-FA-lm and Lanczos-OR-lm.

\section{Numerical experiments and comparison to related algorithms}
\label{sec:comparison}

We now provide several examples which illustrate various aspects of the convergence properties of Lanczos-OR and Lanczos-OR based algorithms, and show when these new methods can outperform more standard techniques like the classic Lanczos-FA.

\subsection{The matrix sign function}
\label{sec:experiment_sign_function}

As we noted in \cref{sec:sign_function}, Lanczos-OR can be used to obtain an approximation to the matrix sign function. 
A related approach, which interpolates the sign function at the so called ``harmonic Ritz values'', is described in \cite[Section 4.3]{eshof_frommer_lippert_schilling_van_der_vorst_02}.
The harmonic Ritz values are characterized by the generalized eigenvalue problem,
$
     [\widehat{\vec{T}}^2]_{:k,:k} \vec{y} = \theta \vec{T} \vec{y},
$
and are closely related to MINRES, which produces a polynomial interpolating $1/x$ at the harmonic Ritz values \cite{paige_parlett_vandervorst_95}.

\begin{example}

We construct a matrix with $400$ eigenvalues, 100 of which are the negatives of the values of a model problem \cite{strakos_91,strakos_greenbaum_92} with parameters $\kappa = 10^2$, $\rho=0.9$, and $n=100$ and 300 of which are the values of a model problem with parameters $\kappa = 10^3$, $\rho=0.8$, $n=300$.
Here, the model problem eigenvalues are given by
\begin{equation*}
    \lambda_1 = 1
    ,\quad \lambda_n = \kappa
    ,\quad \lambda_i = \lambda_1 + \left( \frac{i-1}{n-1} \right) \cdot (\kappa -1) \cdot \rho^{n-i}, \qquad i=2,\ldots,n-1.
\end{equation*}

We compute the Lanczos-OR approximation, the Lanczos-FA approximation, the harmonic Ritz value based approximation from \cite{eshof_frommer_lippert_schilling_van_der_vorst_02}, and the optimal $\vec{A}^2$-norm approximation to the matrix sign function.
The results are shown in \cref{fig:sign_opt}.
In all cases, we use the Lanczos algorithm with full reorthogonalization.
Because eigenvalues of $\vec{T}$ may be near to zero, Lanczos-FA exhibits oscillatory behavior.
On the other hand, the Lanczos-OR based approach and the harmonic Ritz value based approach have much smoother convergence.
Note that the Lanczos-OR induced approximation is not optimal, although it seems to perform close to optimally after a few iterations.

\begin{figure}
    \centering
    \includegraphics[width=.95\textwidth]{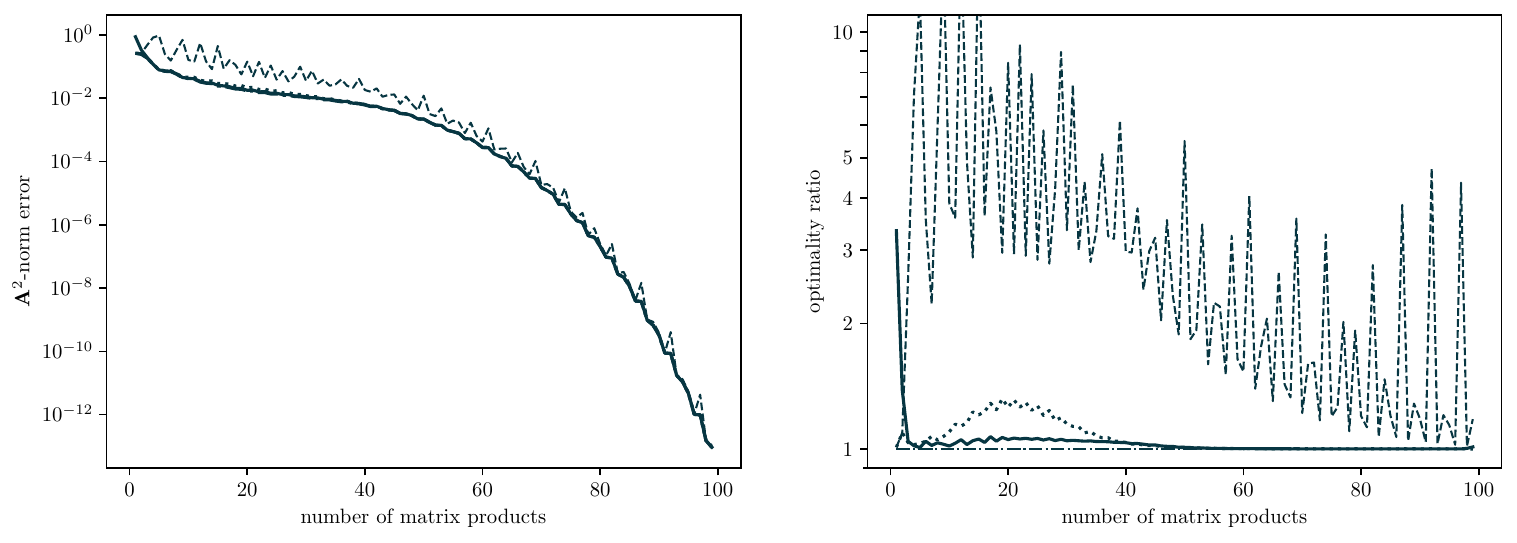}
    \caption{Comparison of $\vec{A}^2$-norm errors for approximating $\operatorname{sign}(\vec{A})\vec{b}$ (normalized by $\|\vec{b}\|_{\vec{A}^2}$).
    \emph{Legend}: Lanczos-OR induced approximation ({\protect\raisebox{0mm}{\protect\includegraphics[scale=.7]{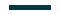}}}),
    interpolation at harmonic Ritz values ({\protect\raisebox{0mm}{\protect\includegraphics[scale=.7]{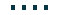}}}),
    Lanczos-FA ({\protect\raisebox{0mm}{\protect\includegraphics[scale=.7]{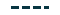}}}),
    $\vec{A}^2$-norm optimal ({\protect\raisebox{0mm}{\protect\includegraphics[scale=.7]{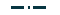}}}).
    \emph{Left}: $\vec{A}^2$-norm of error. 
    \emph{Right}: optimality ratio. 
    \emph{Remark}: Lanczos-OR exhibits smoother convergence than Lanczos-FA and is nearly optimal.
    }
    \label{fig:sign_opt}
\end{figure}

\end{example}

\begin{example}

In this example, we show the spectrum approximations induced by the algorithms from the previous example.
We now set $\vec{A}$ to be a diagonal matrix with $1000$ eigenvalues set to the quantiles of a Chi-squared distribution with parameters $\alpha=1$ and $\beta=10$. 
We set $k=10$ and consider approximations to the function $c\mapsto \vec{b}^\T \bOne [\vec{A} \leq c] \vec{b}$ for a range of values $c$. 
Here $\bOne[x\leq c] = (1-\operatorname{sign}(x-c))/2$ is one if $x\leq c$ and zero otherwise.
We pick $\vec{b}$ as a unit vector with equal projection onto each eigencomponent so that $\vec{b}^\T \bOne [\vec{A} \leq c] \vec{b}$ gives the fraction of eigenvalues of $\vec{A}$ below $c$.
In the $n\to\infty$ limit, this function will converge pointwise to the cumulative distribution of a Chi-squared random distribution with parameters $\alpha=1$ and $\beta=10$. 
The results are shown in \cref{fig:spec_approx}.

Note that the Lanczos-FA based approach is piecewise constant with jumps at each eigenvalue of $\vec{T}$.
On the other hand, the harmonic Ritz value and Lanczos-OR based approaches produce continuous approximations to the spectrum.
In this particular example, the spectrum of $\vec{A}$ is near to a smooth limiting density, so the harmonic Ritz value and Lanczos-OR based approaches seem to produce better approximations.

In general it is not possible to pick $\vec{b}$ with equal projection onto each eigencomponent since the eigenvectors of $\vec{A}$ are unknown. 
However, by choosing $\vec{b}$ from a suitable distribution, it can be guaranteed that $\vec{b}$ has roughly equal projection onto each eigencomponent.
In this case, the Lanczos based approach above is referred to as stochastic Lanczos quadrature \cite{chen_trogdon_ubaru_22}.

\begin{figure}
    \centering
    \includegraphics[width=.95\textwidth]{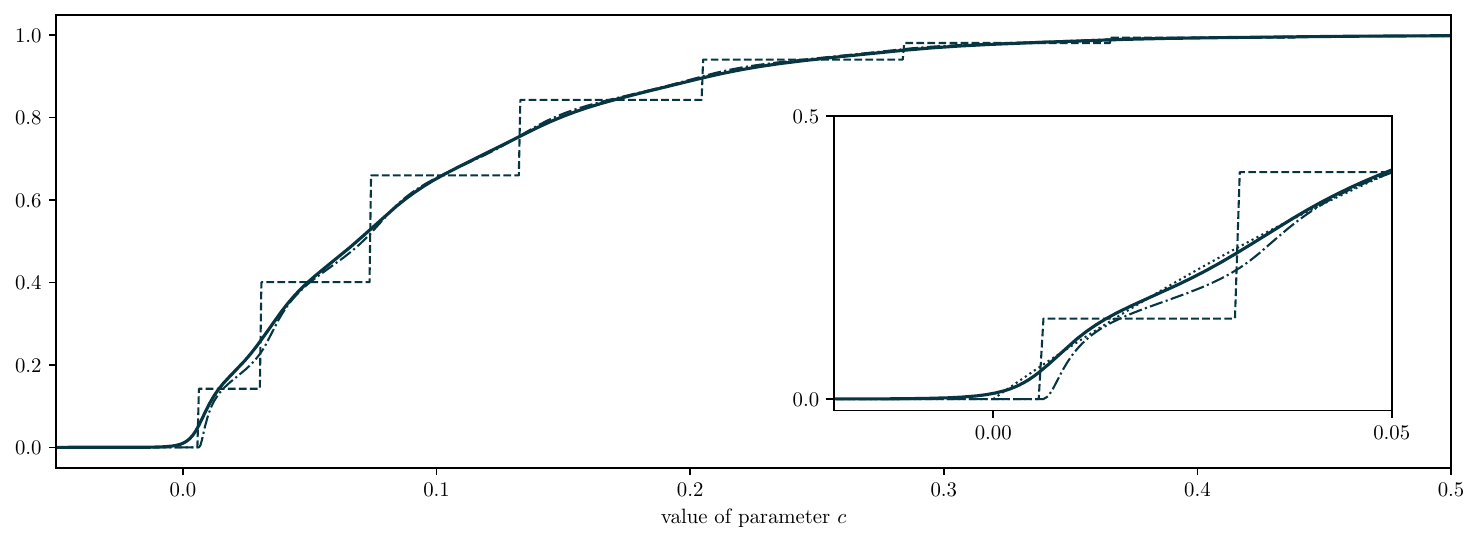}
    \caption{Comparison of spectrum approximations.
        \emph{Legend}: Lanczos-OR induced approximation ({\protect\raisebox{0mm}{\protect\includegraphics[scale=.7]{imgs/legend/blue_solid.pdf}}}),
        Lanczos-FA ({\protect\raisebox{0mm}{\protect\includegraphics[scale=.7]{imgs/legend/blue_dash_thin.pdf}}}),
        harmonic Ritz values based approximation ({\protect\raisebox{0mm}{\protect\includegraphics[scale=.7]{imgs/legend/blue_dashdot_thin.pdf}}}),
        limiting density
        ({\protect\raisebox{0mm}{\protect\includegraphics[scale=.7]{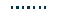}}}).
        .
        \emph{Remark}: The Lanczos-OR and Harmonic Ritz value based approaches provide smooth approximations which match the smooth limiting density much better than the piecewise constant approximation computed with Lanczos-FA.
        }
    \label{fig:spec_approx}
\end{figure}
\end{example}

\subsection{Rational matrix functions}

We now illustrate the effectiveness of the Lanczos-OR based approach to approximating rational matrix functions described in \cref{sec:rational_func}.
Then we compare an existing low-memory approach, called multishift CG, to the analogous approaches based on Lanczos-OR-lm and Lanczos-FA-lm.

Throughout, this section, we will assume that $r$ is a rational function of the form
\begin{equation}
\label{eqn:rat_form}
r(x) = \sum_{i=1}^{m} \frac{A_ix^2+B_ix+C_i}{a_ix^2+b_ix+c_i}.
\end{equation}
This is relatively general since any real valued rational function with numerator degree smaller than denominator degree and only simple poles can be written in this form (in fact, this would be true even if $A_i=0$).
A range of rational functions of this form appear naturally; for instance by a quadrature approximation to a Cauchy integral formula representation of $f$ \cite{hale_higham_trefethen_08}.
Similar rational functions are seen in \cite{eshof_frommer_lippert_schilling_van_der_vorst_02,frommer_simoncini_09}. For rational functions of this form, 
it is clear that $r(\vec{A})\vec{b}$ has the form
\begin{equation}
\label{eqn:rat_form_matrix}
    r(\vec{A})\vec{b} = \sum_{i=1}^{m} (A_i\vec{A}^2 + B_i\vec{A} + C_i\vec{I}) \vec{x}_i 
\end{equation}
where $\vec{x}_i$ is obtained by solving the linear system of equations $(a_i\vec{A}^2 + b_i \vec{A} + c_i \vec{I}) \vec{x}_i = \vec{b}$, and in certain cases, the shift invariance of Krylov subspace can be used to simultaneously compute all of the $\vec{x}_i$ using the same number of matrix-vector products as would be required to approximate a single $\vec{x}_i$ \cite{eshof_frommer_lippert_schilling_van_der_vorst_02,frommer_simoncini_08,frommer_simoncini_08,guttel_schweitzer_21,pleiss_jankowiak_eriksson_damle_gardner_20}.

\begin{example}

In this example, we use the same spectrum as in the first example.
However, rather than approximating the sign function directly, we instead use Lanczos-OR to approximate each term of a proxy rational function of the form \cref{eqn:rat_form}.
In particular, we consider the best uniform approximation\footnote{Note that the eigenvalues of $\vec{A}$ lie in $[-10^2,-1]\cup[1,10^3]$, so we could use an asymmetric approximation to the sign function. 
This would reduce the degree of the rational function required to obtain an approximation of given accuracy, but the qualitative behavior of Lanczos-OR-lm would not change substantially.} of degree $(39,40)$ to the sign function on $[-10^3,-1]\cup[1,10^3]$.
Such an approximation is due to Zolotarev \cite{zolotarev_77}, an can be derived from the more well known Zolotarev approximation to the inverse square root function on $[1,10^6]$.
Our implementation follows the partial fractions implementation in the Rational Krylov Toolbox \cite{berljafa_elsworth_guttel_20} and involves computing the sum of $20$ terms of degree $(1,2)$.
The results are shown in \cref{fig:proxy_rat}.

\begin{figure}[htb]
    \label{fig:proxy_rat}
    \centering
    \includegraphics[width=.95\textwidth]{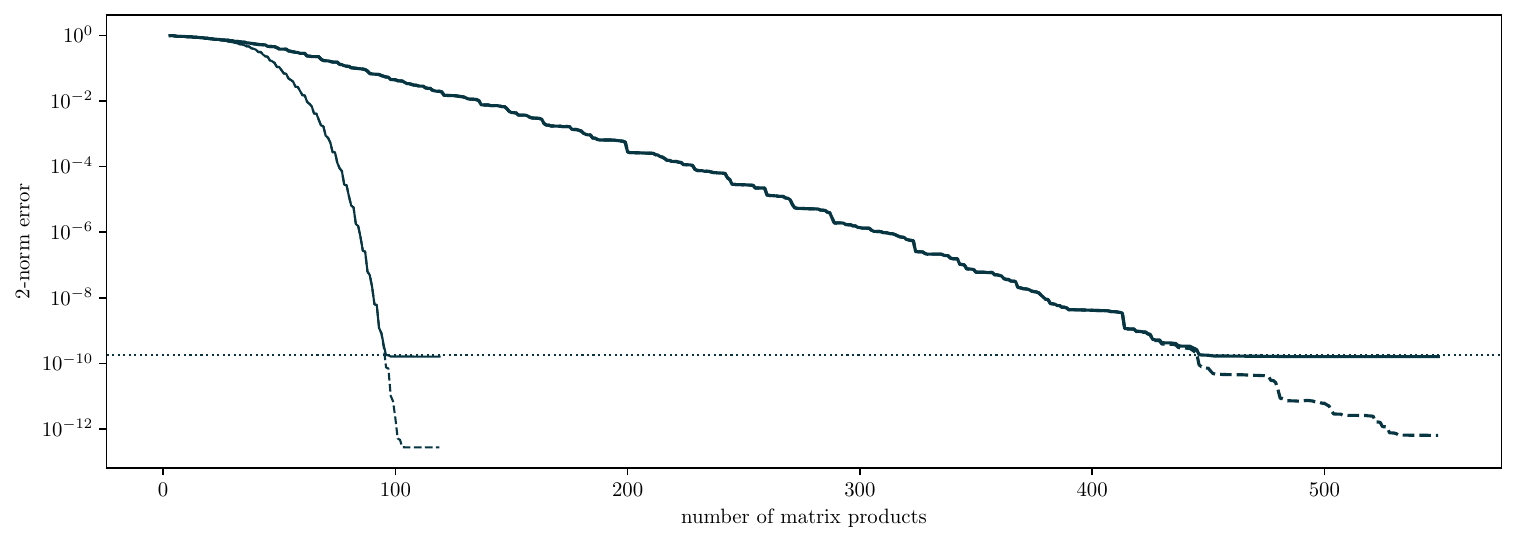}
    \caption{2-norm error in Lanczos-OR-lm based rational approximation ($R(x)=1$) to matrix sign function (normalized by $\|\vec{b}\|_2$).
        \emph{Legend}:
        Lanczos-OR-lm based approximation of matrix sign function without/with reorthogonalization
        ({\protect\raisebox{0mm}{\protect\includegraphics[scale=.7]{imgs/legend/blue_solid.pdf}}}/{\protect\raisebox{0mm}{\protect\includegraphics[scale=.7]{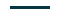}}}),
        Lanczos-OR-lm based approximation of proxy rational matrix function without/with reorthogonalization
        ({\protect\raisebox{0mm}{\protect\includegraphics[scale=.7]{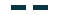}}}/{\protect\raisebox{0mm}{\protect\includegraphics[scale=.7]{imgs/legend/blue_dash_thin.pdf}}}),
        Infinity norm error of proxy rational function approximation ({\protect\raisebox{0mm}{\protect\includegraphics[scale=.7]{imgs/legend/blue_dot_thin.pdf}}}).
        \emph{Remark}:
        The convergence of Lanczos-OR to the sign function matches the convergence to the proxy rational function until the error is of the order of the error in the proxy rational function. 
        }
\end{figure}
\end{example}

At least while the error of the Lanczos-OR approximation to the proxy rational matrix function is large relative to the sign function approximation error, as seen in \cref{eqn:triangle_ineq}, the error in approximating the matrix sign function is similar to the error in approximating the proxy rational matrix function.
However, the final accuracy is limited by the quality of the scalar approximation.
Also note that it really only makes sense to use Lanczos-OR-lm with a short-recurrence version of Lanczos, in which case the effects of a perturbed Lanczos recurrence are prevalent. 
In particular, the algorithm encounters a delay of convergence as compared to what would happen with reorthogonalization.
This is because the example problem's spectrum has many outlying eigenvalues, so the Lanczos algorithm quickly loses orthogonality and begins to find ``ghost eigenvalues'' \cite{meurant_strakos_06,liesen_strakos_13}.

\subsubsection{Comparison of Lanczos-OR, Lanczos-FA, and CG}

To compute terms of \cref{eqn:rat_form_matrix} on could use Lanczos-OR, Lanczos-FA, or assuming the denominator is positive definite, CG (where each CG iteration requires a product with the denominator).
The following example highlights some of the tradeoffs:

\begin{example}
\label{ex:lanczos_msCG_comparison}
We construct several test problems by placing eigenvalues uniformly throughout the specified intervals.
In all cases, $\vec{b}$ has uniform weight onto each eigencomponent.
The outputs are computed using standard Lanczos, but we note that the spectrum and number of iterations are such that the behavior is quite similar to if full reorthgonalization were used. 
In particular, orthogonality is not lost since no Ritz value converges. 
The results of our experiments are shown in \cref{fig:msCG_squared}.

\begin{figure}
    \centering
    \includegraphics[width=.95\textwidth]{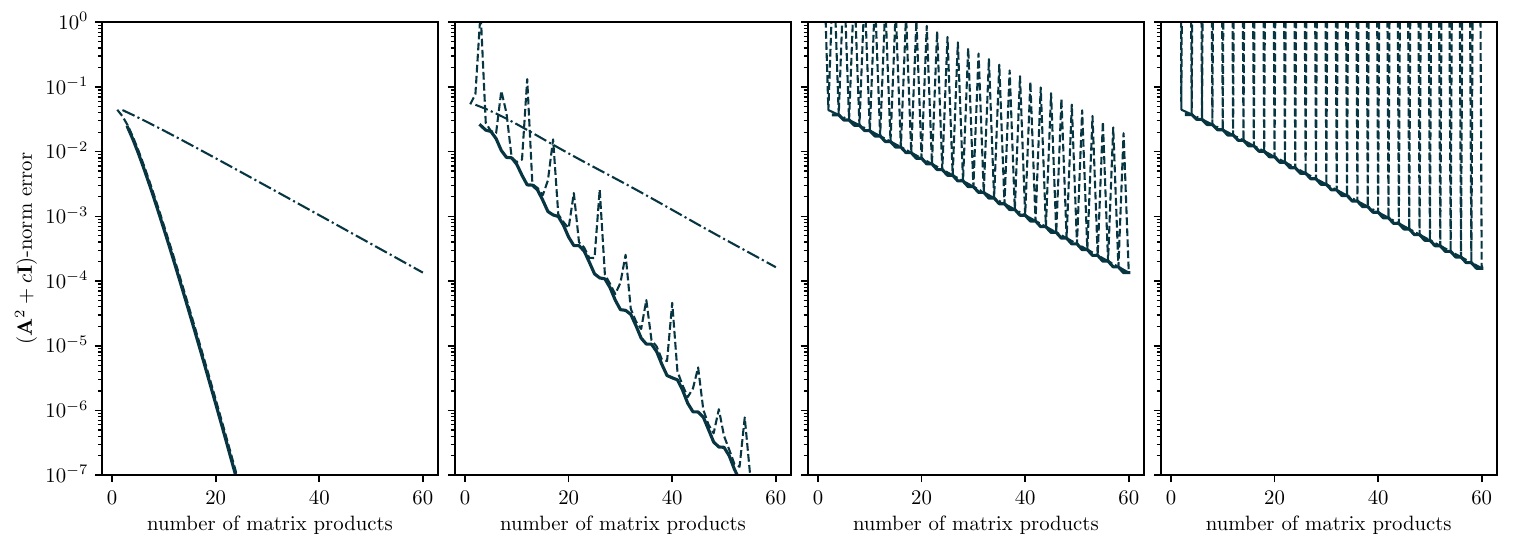}
    \caption{Comparison of $(\vec{A}^2+c\vec{I})$-norm errors for Lanczos-OR ($R(x)=1$), Lanczos-FA, and CG for computing $(\vec{A}^{2}+c\vec{I})^{-1}\vec{b}$ (normalized by $\|\vec{b}\|_{\vec{A}^2+c\vec{I}}$).  
    Here CG works with $\vec{A}^2+c\vec{I}$ and requires two matrix-vector products per iteration whereas Lanczos-FA works with $\vec{A}$ and requires just one.
    \emph{Legend}: Lanczos-OR ({\protect\raisebox{0mm}{\protect\includegraphics[scale=.7]{imgs/legend/blue_solid.pdf}}}),
    Lanczos-FA ({\protect\raisebox{0mm}{\protect\includegraphics[scale=.7]{imgs/legend/blue_dash_thin.pdf}}}),
    CG on squared system ({\protect\raisebox{0mm}{\protect\includegraphics[scale=.7]{imgs/legend/blue_dashdot_thin.pdf}}}).
    \emph{Far left}: eigenvalues on $[1,10]$, $c=0.05$. 
    \emph{Middle left}: eigenvalues on $[-1.5,-1]\cup[1,10]$, $c=0.05$.
    \emph{Middle right}: eigenvalues on $[-10,-1]\cup[1,10]$, $c=0.05$. 
    \emph{Far right}: eigenvalues on $[-10,-1]\cup[1,10]$, $c=0$.
    \emph{Remark}: Lanczos-OR converges without oscillations while automatically matching the rate of convergence of the better of Lanczos-FA and CG on the squared system.
    }
    \label{fig:msCG_squared}
\end{figure}
\end{example}
In the first three examples, we consider approximations to $r(x) = 1/(x^2+0.05)$ with eigenvalues spaced with increments of $0.005$ in $[1,10]$, $[-1.5,-1]\cup[1,10]$, and $[-10,-1]\cup[1,10]$ respectively.
For all these examples, the condition number of $\vec{A}^2+0.05\vec{I}$ is roughly 100 and the eigenvalues of $\vec{A}^2+0.05 \vec{I}$ fill out the interval $[1.05,100.05]$.
As such we observe that multishift CG converges at a rate (in terms of matrix products with $\vec{A}$) of roughly $\exp(-k/\sqrt{\kappa(\vec{A}^2)}) = \exp(-k/\sqrt{100})$ on all of the examples.

In the first example, $\vec{A}$ is positive definite.
Here Lanczos-FA and Lanczos-OR converge similarly to CG on $\vec{A}$ at a rate of roughly $\exp(-2k/\sqrt{10})$, where $k$ is the number of matrix-vector products with $\vec{A}$.

In the next example $\vec{A}$ is indefinite.
The convergence of CG is unchanged, because CG acts on $\vec{A}^2 + c\vec{I}$, it is unable to ``see'' the asymmetry in the eigenvalues of $\vec{A}$.
While the convergence of Lanczos-FA and Lanczos-OR is slowed considerably, both methods converges more quickly than CG due to the asymmetry in the intervals to the left and the right of the origin.
The convergence of these methods is at a rate of roughly $\exp(-k/\sqrt{15})$, although the exact rate is more complicated to compute \cite{fischer_96,schiefermayr_11}.
We also note the emergence of oscillations in the error curve of Lanczos-FA.

In the third example, the asymmetry in the eigenvalue distribution about the origin is removed, and Lanczos-FA and Lanczos-OR converge at a very similar rate to multishift CG.
Note that Lanczos-FA displays larger oscillations, since the symmetry of the eigenvalue distribution of $\vec{A}$ ensures that $\vec{T}$ has an eigenvalue at zero whenever $k$ is odd.
However, the size of the oscillations is regularized by the fact that $c>0$.

In the final example, we use the same eigenvalue distribution as the third example, but now apply the function $r(x) = x^{-2}$.
Here CG and Lanczos-OR behave essentially the same, but the behavior of Lanczos-FA becomes far more oscillatory.
Indeed, the lack of the regularizing term $c\vec{I}$ means that $r(\vec{T}) = \vec{T}^{-2}$ is not even defined when $\vec{T}$ has an eigenvalue at zero.
Lanczos-FA-lm will break down in such settings, as the LDL factorization of $\vec{T}^2$ is not well defined.
Even in less extreme situations, the LDL factorization may become inaccurate.

\subsection{Optimality in the 2-norm}
\label{sec:2normopt}

The Lanczos-OR iterates are optimal in the $\vec{H}$-norm, where $\vec{H} = N(\vec{A})R(\vec{A})$. 
In many situations (including the special cases of CG or MINRES which are respectively optimal in the $\vec{A}$ and $\vec{A}^2$ norms), it is more desirable to have a good approximation in a different norm.
Thus, it is important to understand how the Lanczos-OR iterates behave in other norms, and for concreteness, we focus on the 2-norm.
While the Lanczos-OR iterates cannot be expected to be optimal in the 2-norm, as seen in \cref{thm:2normopt}, they are optimal up to a factor $\sqrt{\kappa(\vec{H})}$.
In many situations, we find that the iterates tend to satisfy a similar bound, but with $\sqrt{\kappa(\vec{H})}$ replaced by some small value (e.g., 2).
However, we believe the $\sqrt{\kappa(\vec{H})}$ factor is necessary in the worst case. 
Thus, for problems where $\vec{H}$ is very poorly conditioned, Lanczos-OR cannot necessarily be guaranteed to output an iterate which is near to the 2-norm optimal iterate.

\begin{figure}[htb]
    \centering
    \includegraphics[width=.95\textwidth]{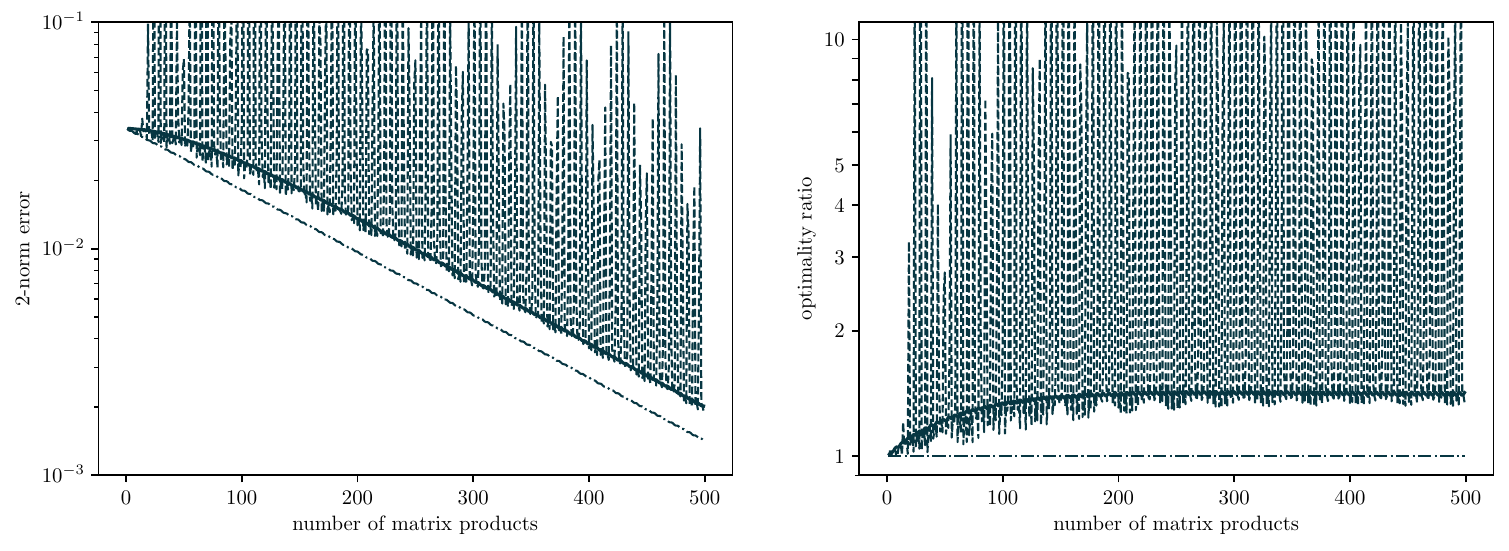}
    \caption{Comparison of Euclidian norm errors for Lanczos-OR and  Lanczos-FA for computing $(\vec{A}^{2}+c\vec{I})^{-1}\vec{b}$ (normalized by $\|\vec{b}\|$) to optimal approximation.
    \emph{Legend}: Lanczos-OR ({\protect\raisebox{0mm}{\protect\includegraphics[scale=.7]{imgs/legend/blue_solid.pdf}}}),
    Lanczos-FA ({\protect\raisebox{0mm}{\protect\includegraphics[scale=.7]{imgs/legend/blue_dash_thin.pdf}}}),
    Optimal Euclidian norm iterate
    ({\protect\raisebox{0mm}{\protect\includegraphics[scale=.7]{imgs/legend/blue_dashdot_thin.pdf}}}).
    \emph{Left}: 2-norm errors. 
    \emph{Right}: optimality ratio.
    \emph{Remark}: Lanczos-OR may perform nearly-optimally, even in the Euclidian norm.
    }
    \label{fig:2normopt}
\end{figure}

\begin{example}
We use a similar setup as in \cref{ex:lanczos_msCG_comparison}.
Specifically we consider the approximation to $(\vec{A}^{2}+0.05\vec{I})^{-1}\vec{b}$, where $\vec{A}$ has $n=109602$ eigenvalues spaced uniformly with spacing $0.005$ in $[-50,-1]\cup[1,500]$.
In \cref{fig:2normopt} we show the 2-norm of the errors for Lanczos-OR and Lanczos-FA in comparison to the optimal 2-norm approximation. 
We also show the optimally ratio, which illustrates that both algorithms perform nearly optimally, although the Lanczos-OR iterates are less erratic.
In particular, the approximation ratio of Lanczos-OR is far smaller than $\sqrt{\kappa(\vec{H})} \approx 500$ for this this particular problem.

\end{example}

\section{Outlook}

There are a range of interesting directions for future work.
A few of the most interesting are summarized here:
\begin{itemize}
    \item In the case $f(x) = 1/x$, \cite{cullum_greenbaum_96} provides a exact relation between CG and MINRES residuals. 
        Can we relate the errors of Lanczos-OR and Lanczos-FA in general? 
    \item Can we provide a sharper comparison between the Lanczos-OR and Lanczos-FA approximations to the matrix sign function?
    \item Is the induced algorithm for the matrix sign function nearly optimal/can we derive simple spectrum dependent bounds?
    \item For what other functions can we use Lanczos-OR to induce a new algorithm?
    \item How does Lanczos-OR generalize ``harmonic Ritz values'', and can this perspective provide any insight into Lanczos-FA?
    \item Can we provide a unified analysis of Krylov subspace methods such as MINRES and CG in finite precision arithmetic?
    \item Why does Lanczos-FA tend to perform ``nearly optimally'', at least in the sense of the smallest error observed at all iterations up to the current iteration?
\end{itemize}

\clearpage
\bibliography{lanczos_rational_opt}

\begin{thebibliography}{10}

\bibitem{afanasjew_eiermann_ernst_guttel_08}
{\sc M.~Afanasjew, M.~Eiermann, O.~G. Ernst, and S.~G\"{u}ttel}, {\em
  Implementation of a restarted krylov subspace method for the evaluation of
  matrix functions}, Linear Algebra and its Applications, 429 (2008),
  pp.~2293--2314.

\bibitem{berljafa_elsworth_guttel_20}
{\sc M.~Berljafa, S.~Elsworth, and S.~Güttel}, {\em A rational {K}rylov
  toolbox for matlab}, 2020.

\bibitem{boricci_00a}
{\sc A.~Bori{\c{c}}i}, {\em Fast methods for computing the {N}euberger
  operator}, Springer Berlin Heidelberg, 2000, pp.~40--47.

\bibitem{chen_greenbaum_musco_musco_22}
{\sc T.~Chen, A.~Greenbaum, C.~Musco, and C.~Musco}, {\em Error bounds for
  lanczos-based matrix function approximation}, {SIAM} Journal on Matrix
  Analysis and Applications, 43 (2022), pp.~787--811.

\bibitem{chen_trogdon_ubaru_22}
{\sc T.~Chen, T.~Trogdon, and S.~Ubaru}, {\em Randomized matrix-free quadrature
  for spectrum and spectral sum approximation}, 2022.

\bibitem{cullum_greenbaum_96}
{\sc J.~Cullum and A.~Greenbaum}, {\em Relations between {G}alerkin and
  norm-minimizing iterative methods for solving linear systems}, SIAM Journal
  on Matrix Analysis and Applications, 17 (1996), pp.~223--247.

\bibitem{druskin_greenbaum_knizhnerman_98}
{\sc V.~Druskin, A.~Greenbaum, and L.~Knizhnerman}, {\em Using nonorthogonal
  {L}anczos vectors in the computation of matrix functions}, SIAM Journal on
  Scientific Computing, 19 (1998), pp.~38--54.

\bibitem{druskin_knizhnerman_89}
{\sc V.~Druskin and L.~Knizhnerman}, {\em Two polynomial methods of calculating
  functions of symmetric matrices}, {USSR} Computational Mathematics and
  Mathematical Physics, 29 (1989), pp.~112--121.

\bibitem{fischer_96}
{\sc B.~Fischer}, {\em Polynomial Based Iteration Methods for Symmetric Linear
  Systems}, Vieweg+Teubner Verlag, 1996.

\bibitem{freund_92}
{\sc R.~W. Freund}, {\em Conjugate gradient-type methods for linear systems
  with complex symmetric coefficient matrices}, {SIAM} Journal on Scientific
  and Statistical Computing, 13 (1992), pp.~425--448.

\bibitem{frommer_guttel_schweitzer_14a}
{\sc A.~Frommer, S.~G\"{u}ttel, and M.~Schweitzer}, {\em Convergence of
  restarted {K}rylov subspace methods for {S}tieltjes functions of matrices},
  {SIAM} Journal on Matrix Analysis and Applications, 35 (2014),
  pp.~1602--1624.

\bibitem{frommer_guttel_schweitzer_14}
\leavevmode\vrule height 2pt depth -1.6pt width 23pt, {\em Efficient and stable
  {A}rnoldi restarts for matrix functions based on quadrature}, {SIAM} Journal
  on Matrix Analysis and Applications, 35 (2014), pp.~661--683.

\bibitem{frommer_schweitzer_15}
{\sc A.~Frommer and M.~Schweitzer}, {\em Error bounds and estimates for
  {K}rylov subspace approximations of {S}tieltjes matrix functions}, {BIT}
  Numerical Mathematics, 56 (2015), pp.~865--892.

\bibitem{frommer_simoncini_08}
{\sc A.~Frommer and V.~Simoncini}, {\em Matrix functions}, in Mathematics in
  Industry, Springer Berlin Heidelberg, 2008, pp.~275--303.

\bibitem{frommer_simoncini_09}
\leavevmode\vrule height 2pt depth -1.6pt width 23pt, {\em Error bounds for
  {L}anczos approximations of rational functions of matrices}, in Numerical
  Validation in Current Hardware Architectures, Berlin, Heidelberg, 2009,
  Springer Berlin Heidelberg, pp.~203--216.

\bibitem{greenbaum_89}
{\sc A.~Greenbaum}, {\em Behavior of slightly perturbed {L}anczos and
  conjugate-gradient recurrences}, Linear Algebra and its Applications, 113
  (1989), pp.~7 -- 63.

\bibitem{greenbaum_97}
\leavevmode\vrule height 2pt depth -1.6pt width 23pt, {\em Iterative Methods
  for Solving Linear Systems}, Society for Industrial and Applied Mathematics,
  Philadelphia, PA, USA, 1997.

\bibitem{greenbaum_druskin_knizhnerman_99}
{\sc A.~Greenbaum, V.~Druskin, and L.~A. Knizhnerman}, {\em On solving
  indefinite symmetric linear systems by means of the {L}anczos method}, Zh.
  Vychisl. Mat. Mat. Fiz., 39 (1999), pp.~371--377.

\bibitem{guttel_schweitzer_21}
{\sc S.~G\"{u}ttel and M.~Schweitzer}, {\em A comparison of limited-memory
  {K}rylov methods for {S}tieltjes functions of {H}ermitian matrices}, {SIAM}
  Journal on Matrix Analysis and Applications, 42 (2021), pp.~83--107.

\bibitem{hale_higham_trefethen_08}
{\sc N.~Hale, N.~J. Higham, and L.~N. Trefethen}, {\em Computing \({A}^\alpha,
  \log({A})\), and related matrix functions by contour integrals}, SIAM Journal
  on Numerical Analysis, 46 (2008), pp.~2505--2523.

\bibitem{hestenes_stiefel_52}
{\sc M.~R. Hestenes and E.~Stiefel}, {\em Methods of conjugate gradients for
  solving linear systems}, vol.~49, NBS Washington, DC, 1952.

\bibitem{higham_02}
{\sc N.~J. Higham}, {\em Accuracy and Stability of Numerical Algorithms},
  Society for Industrial and Applied Mathematics, Jan. 2002.

\bibitem{ilic_turner_simpson_09}
{\sc M.~D. Ilic, I.~W. Turner, and D.~P. Simpson}, {\em A restarted {L}anczos
  approximation to functions of a symmetric matrix}, {IMA} Journal of Numerical
  Analysis, 30 (2009), pp.~1044--1061.

\bibitem{liesen_strakos_13}
{\sc J.~Liesen and Z.~Strako{\v{s}}}, {\em Krylov subspace methods: principles
  and analysis}, Numerical mathematics and scientific computation, Oxford
  University Press, 1st ed~ed., 2013.

\bibitem{lopez_simoncini_06}
{\sc L.~Lopez and V.~Simoncini}, {\em Analysis of projection methods for
  rational function approximation to the matrix exponential}, {SIAM} Journal on
  Numerical Analysis, 44 (2006), pp.~613--635.

\bibitem{meurant_strakos_06}
{\sc G.~Meurant and Z.~Strako{\v{s}}}, {\em The {L}anczos and conjugate
  gradient algorithms in finite precision arithmetic}, Acta Numerica, 15
  (2006), pp.~471--542.

\bibitem{musco_musco_sidford_18}
{\sc C.~Musco, C.~Musco, and A.~Sidford}, {\em Stability of the {L}anczos
  method for matrix function approximation}, in Proceedings of the Twenty-Ninth
  Annual ACM-SIAM Symposium on Discrete Algorithms, SODA ’18, USA, 2018,
  Society for Industrial and Applied Mathematics, p.~1605–1624.

\bibitem{niehoff_06}
{\sc J.~Niehoff}, {\em Projektionsverfahren zur Approximation von
  Matrixfunktionen mit Anwendungen auf die Implementierung exponentieller
  Integratoren}, PhD thesis, Heinrich-Heine Universität Düsseldorf,
  Mathematisches Institut, 2006.

\bibitem{paige_71}
{\sc C.~C. Paige}, {\em The computation of eigenvalues and eigenvectors of very
  large sparse matrices.}, PhD thesis, University of London, 1971.

\bibitem{paige_76}
\leavevmode\vrule height 2pt depth -1.6pt width 23pt, {\em {Error Analysis of
  the {L}anczos Algorithm for Tridiagonalizing a Symmetric Matrix}}, IMA
  Journal of Applied Mathematics, 18 (1976), pp.~341--349.

\bibitem{paige_80}
\leavevmode\vrule height 2pt depth -1.6pt width 23pt, {\em Accuracy and
  effectiveness of the {L}anczos algorithm for the symmetric eigenproblem},
  Linear Algebra and its Applications, 34 (1980), pp.~235 -- 258.

\bibitem{paige_parlett_vandervorst_95}
{\sc C.~C. Paige, B.~N. Parlett, and H.~A. Van~der Vorst}, {\em Approximate
  solutions and eigenvalue bounds from krylov subspaces}, Numerical linear
  algebra with applications, 2 (1995), pp.~115--133.

\bibitem{paige_saunders_75}
{\sc C.~C. Paige and M.~A. Saunders}, {\em Solution of sparse indefinite
  systems of linear equations}, {SIAM} Journal on Numerical Analysis, 12
  (1975), pp.~617--629.

\bibitem{pleiss_jankowiak_eriksson_damle_gardner_20}
{\sc G.~Pleiss, M.~Jankowiak, D.~Eriksson, A.~Damle, and J.~R. Gardner}, {\em
  Fast matrix square roots with applications to {G}aussian processes and
  bayesian optimization}, 2020.

\bibitem{saad_92}
{\sc Y.~Saad}, {\em Analysis of some {K}rylov subspace approximations to the
  matrix exponential operator}, SIAM Journal on Numerical Analysis, 29 (1992),
  pp.~209--228.

\bibitem{saad_03}
{\sc Y.~Saad}, {\em Iterative Methods for Sparse Linear Systems}, Society for
  Industrial and Applied Mathematics, Jan. 2003.

\bibitem{schiefermayr_11}
{\sc K.~Schiefermayr}, {\em Estimates for the asymptotic convergence factor of
  two intervals}, Journal of Computational and Applied Mathematics, 236 (2011),
  p.~28–38.

\bibitem{strakos_91}
{\sc Z.~Strakos}, {\em On the real convergence rate of the conjugate gradient
  method}, Linear Algebra and its Applications, 154-156 (1991), pp.~535 -- 549.

\bibitem{strakos_greenbaum_92}
{\sc Z.~Strakos and A.~Greenbaum}, {\em Open questions in the convergence
  analysis of the {L}anczos process for the real symmetric eigenvalue problem},
  University of Minnesota, 1992.

\bibitem{eshof_frommer_lippert_schilling_van_der_vorst_02}
{\sc J.~van~den Eshof, A.~Frommer, T.~Lippert, K.~Schilling, and H.~van~der
  Vorst}, {\em Numerical methods for the {QCD}d overlap operator. {I}.
  sign-function and error bounds}, Computer Physics Communications, 146 (2002),
  pp.~203 -- 224.

\bibitem{zolotarev_77}
{\sc E.~Zolotarev}, {\em Application of elliptic functions to questions of
  functions deviating least and most from zero}, Zap. Imp. Akad. Nauk. St.
  Petersburg, 30 (1877), pp.~1--59.

\bibitem{simonova_tichy_21}
{\sc D.~Šimonová and P.~Tichý}, {\em When does the {L}anczos algorithm
  compute exactly?}, 2021.

\end{thebibliography}
\bibliographystyle{siam}

\end{document}